\documentclass[12pt]{article}%

\usepackage{amsmath}
\usepackage{amsfonts}
\usepackage{amssymb}
\usepackage{graphicx}%
\numberwithin{equation}{section}
\begin{document}

\title{Existence of Metrics with Prescribed Scalar Curvature
on the Volume Element Preserving Deformation}
\author{Zhi-Zhang Wang\thanks{School of Mathematical Sciences, Fudan University, Shanghai,China(200433)}}

\date{}
\maketitle

\begin{center}
{\textbf{Abstract}}
\end{center}
\quad\ In this paper,we obtain two results on closed Reimainnian
manifold $M\times [0,T]$.When $T$ is small enough,to any prescribed
scalar curvature, the existence and uniqueness of metrics are
obtained on the volume element preserving deformation.When $T$ is
large and the given scalar curvature is small enough,the same result
holds.

\section{Introduction}
 \quad\ This article concerns the following question:Suppose $(M,g),(N,h)$
are m-dimension and n-dimension Riemainnian manifolds,on which $g,h$
are corresponded metrics.Then $g+h$ is a metric on $M\times
N$.(Here,
$(g+h)(X_P,Y_P)=g((X_1)_P,(Y_1)_P)+h((X_2)_P,(Y_2)_P),X_P=(X_1)_P+(X_2)_P,Y_P=(Y_1)_P+(Y_2)_P,$
$ X_P,Y_P$ are tangent vectors on $M\times N$ at point $P$.They can
be decomposed into $(X_1)_P,(Y_1)_P$ and $(X_2)_P,(Y_2)_P$, which
are tangent vectors on $M$ and $N$.)The author considers a class of
metric deformation: in fact, one can construct a metric $K$ on
$M\times N$,such that $K=\rho^{n}g+\rho^{-m}h$,here $\rho$ is a
smooth enough positive function defined on $M\times N$. Obvious,
$det(K)=det(\rho^ng)det(\rho^{-m}h)=det(g)det(h)=det(g+h)$, then the
volume forms of $g+h$ and $K$ are same at every point.So we call
this deformation as the volume element preserving deformation. And
the question is that with what condition, a smooth enough function
$\tilde{R}$ defined on $M\times N$ can be viewed as the scalar
curvature of this kind of metrics $K$.

 This question can induce a partial differential
equation namely,
\begin{eqnarray}
\tilde{R}&=&\rho^{-n}R_g+\rho^{m}R_h+\rho^{-n}
\triangle_gln\rho^n+\rho^m\triangle_hln\rho^{-m}\nonumber\\&&
-\frac{n}{4}(nm+2n+m^2)\rho^{-n}|\nabla_gln\rho|_g^2\nonumber\\&&
-\frac{m}{4}(nm+2m+n^2)\rho^m|\nabla_hln\rho|_h^2
\end{eqnarray}
It is too complicated,the author only study this equation when
$n=1$.Do a transformation, this equation becomes,
\begin{eqnarray}
\tilde{R}&=& e^{-2u}R_g+2e^{-2u}\triangle
u-2me^{2mu}\frac{\partial^2 u}{\partial t^2}\nonumber\\&&
-(m^2+m+2)e^{-2u}|\nabla u|^2-m(3m+1)e^{2mu}u^2_t
\end{eqnarray}
It is a degenerate quasi-linear hyperbolic equation.

 This article is composed by 5 sections.In section 2,the
author introduces the notations of some spaces and corresponded
norms, refers some Lemmas which are useful to the rest part of
this paper.Section 3 gives the details to obtain equation
(1.1),(1.2) from the above proposed question.In section 4,using
the energy estimate,the prior estimate of equation(1.2) obtained
on $M\times I$,here $I=[0,T]$ is a interval of length $T$.Section
5 studies the linear equation of the following form:
\begin{equation}
\frac{\partial^2u}{\partial t^2}+a(x,t)\frac{\partial u}{\partial
t}-(\alpha (x,t)\triangle u+<\nabla\beta(x,t),\nabla
u>+\gamma(x,t)u)=f(x,t)
\end{equation}
Using the estimate of section 4 and Banach contraction principle,
the last section gives two results.

In a word,the paper mainly obtains the following two theorems.\\

\noindent\textbf{Theorem 1.1}\quad\ Suppose $(M,g)$ be m-dimensional
closed Riemainnian manifold, $I=[0,T]$ be a interval of length $T$
with standard Euclidean metric $h$, and $\tilde{R}\in
C^{\infty}(M\times I)$.To any integer $k>\frac{m}{2}+3$,exist $0<
t_0\leq T$, then exist a unique volume element preserving
deformation $e^{2u}g+e^{-2mu}h$ such that the scalar curvature of
this metric is $\tilde{R}$ on $M\times [0,t_0]$.Here, $u\in
C^k(M\times [0,t_0])$,and satisfies
$u(x,0)=\varphi(x),u_t(x,0)=\psi(x)$ and $\varphi(x),\psi(x)\in
C^{\infty}(M)$ are two prescribed functions.And $t_0$ relates to
manifold $M,k,\tilde{R},\varphi,\psi$.\\

\noindent\textbf{Theorem 1.2}\quad\  Suppose $(M,g)$ be
m-dimensional closed Riemainnian manifold which is scalar curvature
flat $R_g=0$.And $I=[0,T]$ be a interval of length $T$ with standard
Euclidean metric $h$, and $\tilde{R}\in C^{\infty}(M\times I)$.To
any integer $k>\frac{m}{2}+3$,exist $\varepsilon
>0$,and when
$\|\tilde{R}\|_{H^{[\frac{m}{2}+k]}(M)}\leq \varepsilon$ then exist
a unique volume element preserving deformation such that the scalar
curvature of the metric $e^{2u}g+e^{-2mu}h$ is $\tilde{R}$ on
$M\times I$.Here, $u\in C^k(M\times I)$,and $u$ satisfies
$u(x,0)=0,u_t(x,0)=0$ ,and $\varepsilon$ only relates to manifold
$M,k,I$.\\

\noindent\textbf{Remark 1.3}\quad\ Here,closed manifold means
compact
oriented and without boundary manifold.\\

 At last ,the author should thank to his supervisor
Professor Huang Xuanguo.This question is proposed originally by
him.And he induced the correspond equation in some special cases.The
section 3 is the generalize of his computation.He carefully checks
my paper,and gives me much valuable suggestion.This paper can not be
completed without his continual encouragement.And the author also
wants to thank to his classmates:Du Yi,Hu Junqi,Yang Yongfu.They
give me a lot of help during the paper finished.
\section{Preliminary}
\quad\ Suppose $(M,g)$ be a closed m-dimensional Riemainnian
manifold. $\{ x_i\}^m_{i=1}$ is its local coordinate.$(g^{ij})$ is
the inverse matrix of $(g_{ij})$. Let $u$ be a smooth enough
function on $M$. Using
$\nabla_{i_s}\nabla_{i_{s-1}}\cdots\nabla_{i_1}u$ as $u$'s the
$s$-order covariant derivation,here $\nabla$ the Levi-Civita
connection to $g$.We use the following notations:
\begin{eqnarray}
&&|\nabla^su|^2_g=\sum^m_{i_1,\cdots, i_s;j_1,\cdots,
j_s=1}(g^{i_1j_1}\cdots g^{i_sj_s})(\nabla_{i_s}\cdots
\nabla_{i_1}u)(\nabla_{j_s}\cdots
\nabla_{j_1}u);\nonumber\\&&\|u\|_{L^p(M)}=(\int_M|u|^pdV_g)^{1/p};
\|\nabla^su\|_{L^p(M)}=(\int_M|\nabla^su|^p_gdV_g)^{1/p};\nonumber\\&&
\|u\|_{W^{s,p}(M)}=(\sum^s_{j=0}\|\nabla^ju\|^p_{L^p(M)})^{1/p};\|u\|_{C^0(M)}=max_M|u|;\nonumber\\&&
\|\nabla^ku\|_{C^0(M)}=max_M|\nabla^ku|_g;
\|u\|_{C^k(M)}=(\sum^k_{j=0}\|\nabla^ju\|^2_{C^0(M)})^{1/2}
\end{eqnarray}
here,$1\leq p,s,k<\infty$.Denote $W^{s,p}(M),C^k(M)$ as the Banach
space with the corresponded norms.Especially,denote $H^s(M)$ as
$W^{s,2}(M)$.

Let $I=[0,T]$ be a interval of length $T$,with the metric $h$ and
the arch length coordinate $\{t\}$.Then on $M\times I$ one can
choose a local coordinate $\{x_i\}^{m+1}_{i=1}$ where
$x_{m+1}=t$.Then direct computation shows
$\Gamma^C_{AB}=0,R^D_{ABC}=0$,if some of $A,B,C,D\in \{1,\cdots,
m+1\}$ are $m+1$.Here, $\Gamma,R$ are the christoffel symbol and
curvature operator of product metric $\bar{g}:=g+h$.Using this
result and Ricci identity, one can finds,to any smooth enough
$(r,s)$ tensor field $A$ on $M\times I$ which has the form in local
$A=\sum^m_{i_1,\cdots,i_r;j_1,\cdots,j_s=1}A^{i_1\cdots
i_r}_{j_1\cdots j_s}(x,t)\frac{\partial}{\partial
x^{i_1}}\otimes\cdots\otimes\frac{\partial}{\partial x^{i_r}}\otimes
dx^{j_1}\otimes\cdots\otimes dx^{j_s}$,
 $\bar{\nabla}_iA^{i_1\cdots
i_r}_{j_1\cdots j_s}$ $=\nabla_iA^{i_1\cdots i_r}_{j_1\cdots
j_s},\bar{\nabla}_{\frac{\partial}{\partial t}}A^{i_1\cdots
i_r}_{j_1\cdots j_s}=\frac{\partial }{\partial t}A^{i_1\cdots
i_r}_{j_1\cdots j_s},i\in\{1,\cdots,m\}$,and
\begin{eqnarray}
\frac{\partial}{\partial t}\nabla^lA^{i_1\cdots i_r}_{j_1\cdots
j_s}=\nabla^l\frac{\partial}{\partial t}A^{i_1\cdots i_r}_{j_1\cdots
j_s}\end{eqnarray} here,$\bar{\nabla}$ is the corresponded
connection of $\bar{g}$.Then we introduce the following norms:
\begin{eqnarray}
&\|\partial^k_tu\|_{L^p([0,t],W^{s,q}(M))}=\|\|\bar{\nabla}^k_{\frac{\partial}{\partial
t}}u\|_{W^{s,q}(M)}\|_{L^p[0,t]}=\|\|\frac{\partial^ku}{\partial
t^k}\|_{W^{s,q}(M)}\|_{L^p[0,t]}&\nonumber\\&
\|u\|_{W^{k,p}([0,t],W^{s,q}(M))}=(\sum^k_{i=0}\|\partial^i_tu\|^p_{L^p([0,t],W^{s,q}(M))})^{1/p}&\nonumber\\&\|u\|_{C^k([0,t],C^s(M))}
=(\sum^k_{i=0}\|\|\partial^i_tu\|^2_{C^s(M)}\|_{C^0[0,t]})^{1/2}&
\end{eqnarray}
and the space $L^p([0,t],W^{s,q}(M))$ is the completion of the
$C^{\infty}(M\times [0,t])$ with this norm,$0\leq k,s<\infty,1\leq
p\leq \infty,1\leq q< \infty$, here $u$ is a smooth enough function
on $M\times I$ and $0<t\leq T$. It is weakly star compact.

The following three Lemmas are classical in partial differential
equation,we refer [H\"{o}rmander] p106-108 [Nirenberg] and [Zheng]
p10-11,p186-187 for the detail
proof.\\

\noindent\textbf{Lemma 2.1}(Galiardo-Nirenberg Inequality)\quad\ Let
$j,n\in\mathbb{Z}$ and $0\leq j<n$. Let $1\leq q,r\leq
+\infty,p\in\mathbb{R},\frac{j}{n}\leq a\leq1$,such that
$\frac{1}{p}-\frac{j}{m}=a(\frac{1}{r}-\frac{n}{m})+(1-a)\frac{1}{q}$.For
any $u\in W^{n,r}(\mathbb{R}^m)\bigcap L^q(\mathbb{R}^m)$,there is a
positive constant $C$ depending on $n,m,j,q,r,a$ such that
\begin{eqnarray}
\|\nabla^ju\|_{L^p(\mathbb{R}^m)}&\leq&C\|\nabla^nu\|^a_{L^r(\mathbb{R}^m)}\|u\|^{1-a}_{L^q(\mathbb{R}^m)}
\end{eqnarray}
with a exception:if $1<r<\infty ,n-j-\frac{m}{r}$ is a nonnegative
integer then the inequality hold only $\frac{j}{n}\leq a<1$.\\

\noindent\textbf{Lemma 2.2}\quad\ suppose
$\frac{1}{r}=\frac{1}{p}+\frac{1}{q},1\leq r,p,q\leq \infty$,and
suppose that all norms appearing what follows are bounded.Then for
any integer $s\geq 0$,one has
\begin{eqnarray}
\|\nabla^s(fg)\|_{L^r(\mathbb{R}^m)}&\leq&C(\|f\|_{L^p(\mathbb{R}^m)}\|\nabla^sg\|_{L^q(\mathbb{R}^m)}
+\|\nabla^sf\|_{L^q(\mathbb{R}^m)}\|g\|_{L^p(\mathbb{R}^m)})\nonumber\\
\end{eqnarray}
here $C$ is not depend on $f,g$.\\

\noindent\textbf{Lemma 2.3}\quad\ Suppose $F:\mathbb{R}^1\longmapsto
\mathbb{R}^1$ is a smooth function satisfies $F(0)=0$.For any
integer $s\geq0$,if for $w\in W^{s,p}(\mathbb{R}^m),1\leq
p\leq\infty$ and $\|w\|_{L^{\infty}(\mathbb{R}^m)}\leq\nu_0$,then
$F(w)\in W^{s,p}(\mathbb{R}^m)$ and
\begin{eqnarray}
\|F(w)\|_{W^{s,p}(\mathbb{R}^m)}&\leq&C\|F\|_{C^s([-\nu_0,\nu_0])}(1+\|w\|^s_{L^{\infty}(\mathbb{R}^m)})\|w\|_{W^{s,p}(\mathbb{R}^m)}
\end{eqnarray}
here,$C$ is not depend on $F,w$.\\

\noindent\textbf{Corollary 2.4}\quad\ On a closed Reimainnian
manifold $(M,g)$,replace $\mathbb{R}^m$ by $M$ in the hypothesis of
Lemma 2.1-2.3, then the corresponded result becomes:
\begin{eqnarray}
&(1)&\|\nabla^ju\|_{L^p(M)}\leq C\|u\|^a_{W^{n,r}(M)}\|u\|^{1-a}_{L^q(M)}\nonumber\\
&(2)&\|\nabla^s(fg)\|_{L^r(M)}\leq
C(\|f\|_{L^p(M)}\|g\|_{W^{s,q}(M)}
+\|f\|_{W^{s,q}(M)}\|g\|_{L^p(M)})\nonumber\\&& \|\nabla^s<\nabla
f,\nabla g>_g\|_{L^r(M)}\leq
C(\|f\|_{W^{1,p}(M)}\|g\|_{W^{s+1,q}(M)}\nonumber\\&&
+\|f\|_{W^{s+1,q}(M)}\|g\|_{W^{1,p}(M)})\nonumber\\
&(3)&\|F(w)\|_{W^{s,p}(M)}\leq
C\|F\|_{C^s([-\nu_0,\nu_0])}(1+\|w\|^s_{C^0(M)})\|w\|_{W^{s,p}(M)}
\end{eqnarray}
here $C$ also relates manifold $M$.\\
\textbf{Convention 2.5}\quad\ Here, $W^{k,\infty}(M)$ which appears
in the above understood as $C^k(M)$,namely
$\|\cdot\|_{W^{k,\infty}(M)}=\|\cdot\|_{C^k(M)}$.\\
\textbf{Proof of Corollary 2.4}\quad\ The proof of (1)(2) are
similar,we only prove (1).Since $M$ is compact,$M$ has finite
partition of unit,namely $\{(U_i,\varphi_i)\}^N_{i=1},$
$\sum^{N}_{i=1}\varphi_i=1,\varphi\geq0,\varphi_i\in
C^{\infty}(M),supp\varphi_i\subset U_i$.Since
\begin{eqnarray}
\|\varphi_iu\|_{W^{n,r}(U_i)}&\leq&C(\sum^n_{l=0}\sum^l_{k=0}\|\nabla^ku\|^r_{L^r(U_i)})^{1/r}\nonumber\\&\leq&
C(\sum^n_{k=0}\|\nabla^ku\|^r_{L^r(U_i)})^{1/r}\nonumber\\&\leq&C\|u\|_{W^{n,r}(M)}
\end{eqnarray}
by Lemma 2.1 and the norm equivalence on $U_i$,namely
$\|\cdot\|_{W^{s,p}(U_i)}$ with the metric $g$ equivalence to
$\|\cdot\|_{W^{s,p}(U_i)}$ with the standard Euclidean metric.One
has
\begin{eqnarray}
\|\nabla^j(\varphi_iu)\|_{L^p(U_i)}&\leq&C\|\varphi_iu\|^a_{W^{n,r}(U_i)}\|\varphi_iu\|^{1-a}_{L^q(U_i)}
\end{eqnarray}
therefoce,using (2.8)(2.9),
\begin{eqnarray}
\|\nabla^ju\|_{L^p(M)}&\leq&\sum^N_{i=1}\|\nabla^j(\varphi_iu)\|_{L^p(U_i)}\nonumber\\
&\leq&C\sum^N_{i=1}\|u\|^a_{W^{n,r}(M)}\|u\|^{1-a}_{L^q(M)}\nonumber\\&\leq&C\|u\|^a_{W^{n,r}(M)}\|u\|^{1-a}_{L^q(M)}
\end{eqnarray}
This gives the result.To (3) the proof is similar to Lemma 2.3 since
(1) is right.\\

\noindent\textbf{Remark 2.6}\quad\ If the hypothesis of Corollary
2.4(1) becomes $1\leq j\leq n,1\leq
q,r\leq\infty,p\in\mathbb{R},\frac{j-1}{n-1}\leq a\leq1$,
$\frac{1}{p}-\frac{j-1}{m}=a(\frac{1}{r}-\frac{n-1}{m})+(1-a)\frac{1}{q}$,
then (2.7)(1) becomes, \begin{eqnarray} \|\nabla^ju\|_{L^p(M)}\leq
C\|u\|^a_{W^{n,r}(M)}\|u\|^{1-a}_{W^{1,q}(M)}
\end{eqnarray} and with the
similar exception:if $1<r<\infty ,n-j-\frac{m}{r}$ is a nonnegative
integer then the inequality hold
only $\frac{j-1}{n-1}\leq a<1$.\\

At last,we refer the following Lemma which is need in section 5.The
detail proof can be found in [Zheng] chapter 3 p103.\\

\noindent\textbf{Lemma 2.7}\quad\ Let $B,B^*$ be both Banach space,
and $B$ is the dual of $B^*$.Suppose $1< p\leq \infty$ and
$u_n\longrightarrow u$ weakly star in
$L^p([0,t],B)$,$u^{\prime}_n\longrightarrow u^{\prime}$ weakly star
in $L^p([0,t],B)$,then $u_n(0)\longrightarrow u(0)$ weakly star in
$B$.Here,$u^{\prime}_n,u^{\prime}$ are the derivation in the meaning
of distribution(or current).
\section{Induced Equation} \quad\ The hypothesis of $(M,g),(N,h),\rho$ as in section
1, $\Sigma=M\times N$,then $g+h$ is the product metric on
$\Sigma$.Then one can construct a new metric
$K=\rho^ng+\rho^{-m}h$.To some point $P$,let $\{x_i\}^m_{i=1}$ and
$\{y_{\alpha}\}^{m+n}_{\alpha=m+1}$ be the normal coordinate on
$(M,g)$ and $(N,h)$,we conserve that $i,j,k,\cdots\in\{1\cdots
m\},\alpha,\beta,\gamma,\cdots\in\{m+1\cdots m+n\}$ and
$A,B,C,\cdots\in\{1,\cdots,m+n\}$.So,
\begin{eqnarray}
g_{ij}(P)=\delta_{ij};\Gamma^k_{ij}(P)=0;h_{\alpha\beta}(P)=\delta_{\alpha\beta};
\Gamma^{\gamma}_{\alpha\beta}(P)=0
\end{eqnarray}
with this coordinate,
\begin{eqnarray}
K_{ij}=\rho^ng_{ij};K_{\alpha\beta}=\rho^{-m}h_{\alpha\beta};K_{i\alpha}=0
\end{eqnarray}
 Denote $(g^{ij}),(h^{\alpha\beta}),(K^{AB})$ the inverse matrix of
$(g_{ij}),(h_{\alpha\beta}),(K_{AB})$ respectively.And denote
$\tilde{\Gamma}^A_{BC},\Gamma^k_{ij},\Gamma^{\gamma}_{\alpha\beta}$
the christoffel symbols to $K,g,h$ respectively. At point $P$,
\begin{eqnarray}
\tilde{\Gamma}_{ABC}&=&\frac{1}{2}(\frac{\partial K_{AC}}{\partial
x_B}+\frac{\partial K_{BC}}{\partial x_A}-\frac{\partial
K_{AB}}{\partial x_C})
\\\tilde{\Gamma}_{ijk}&=&\frac{1}{2}(\frac{\partial
K_{ik}}{\partial x_j}+\frac{\partial K_{jk}}{\partial
x_i}-\frac{\partial K_{ij}}{\partial
x_k})\nonumber\\&=&\rho^n\Gamma_{ijk}+ \frac{1}{2}(\frac{\partial
\rho^n}{\partial x_j}g_{ik}+\frac{\partial \rho^n}{\partial
x_i}g_{jk}-\frac{\partial \rho^n}{\partial x_k}g_{ij})
\end{eqnarray}
then,at point $P$,
\begin{eqnarray}
\tilde{\Gamma}^k_{ij}&=&\sum^m_{l=1}K^{kl}\tilde{\Gamma}_{ijl}=\rho^{-n}\sum^m_{l=1}g^{kl}
\tilde{\Gamma}_{ijl}\nonumber\\&=&\Gamma^k_{ij}+\frac{1}{2}(\frac{\partial
ln\rho^n}{\partial x_j}\delta_{ik}+\frac{\partial ln\rho^n}{\partial
x_i}\delta_{jk}-\sum^m_{l=1}\frac{\partial ln\rho^n}{\partial
x_l}g_{ij}g^{kl})\\
\tilde{\Gamma}_{ij\alpha}&=&\frac{1}{2}(-\frac{\partial
K_{ij}}{\partial y_{\alpha}})=-\frac{1}{2}\frac{\partial
\rho^n}{\partial
y_{\alpha}}g_{ij}\\\tilde{\Gamma}^{\alpha}_{ij}&=&\sum^{m+n}_{\beta=m+1}K^{\alpha\beta}
\tilde{\Gamma}_{ij\beta}=\rho^m\sum^{m+n}_{\beta=m+1}h^{\alpha\beta}(-\frac{1}{2}\frac{\partial
\rho^n}{\partial y_{\beta}})g_{ij}\\\tilde{\Gamma}_{\alpha
jk}&=&\frac{1}{2}\frac{\partial K_{jk}}{\partial
y_{\alpha}}=\frac{1}{2}\frac{\partial\rho^n}{\partial
y_{\alpha}}g_{jk} \\ \tilde{\Gamma}^k_{\alpha
j}&=&\sum^m_{l=1}K^{kl}\tilde{\Gamma}_{\alpha
jl}=\frac{1}{2}\frac{\partial ln\rho^n}{\partial
y_{\alpha}}\delta_{jk}\\\tilde{\Gamma}_{i\beta\gamma}&=&\frac{1}{2}\frac{\partial
K_{\beta\gamma}}{\partial x_i}=\frac{1}{2}\frac{\partial
\rho^{-m}}{\partial
x_i}h_{\beta\gamma}\end{eqnarray}\begin{eqnarray}\tilde{\Gamma}^{\gamma}_{i\beta}&=&\sum^{m+n}_{\delta=m+1}K^{\gamma\delta}
\tilde{\Gamma}_{i\beta\delta}=\frac{1}{2}\frac{\partial
ln\rho^{-m}}{\partial
x_i}\delta_{\gamma\beta}\\\tilde{\Gamma}^{\gamma}_{\alpha\beta}&=&\Gamma^{\gamma}_{\alpha\beta}
+\frac{1}{2}(\frac{\partial ln\rho^{-m}}{\partial
y_{\beta}}\delta_{\alpha\gamma}+\frac{\partial ln\rho^{-m}}{\partial
y_{\alpha}}\delta_{\beta\gamma}-\sum^{m+n}_{\delta=m+1}\frac{\partial
ln\rho^{-m}}{\partial
y_{\delta}}h_{\alpha\beta}h^{\gamma\delta})\nonumber\\\\\tilde{\Gamma}_{\alpha\beta
k}&=&\frac{1}{2}(-\frac{\partial K_{\alpha\beta}}{\partial
x_k})=-\frac{1}{2}\frac{\partial \rho^{-m}}{\partial
x_k}h_{\alpha\beta}\\\tilde{\Gamma}^k_{\alpha\beta}&=&\sum^m_{l=1}K^{kl}\Gamma_{\alpha\beta
l}=\rho^{-n}\sum^m_{l=1}g^{kl}(-\frac{1}{2}\frac{\partial
\rho^{-m}}{\partial x_l})h_{\alpha\beta}
\end{eqnarray}
Using (3.5) and (3.11),at point $P$,
\begin{eqnarray}
\sum^{m+n}_{A=1}\tilde{\Gamma}^A_{iA}&=&\sum^m_{k=1}\tilde{\Gamma}^k_{ik}+\sum^{m+n}_{\beta=m+1}
\tilde{\Gamma}^{\beta}_{i\beta}\nonumber\\&=&\sum^m_{k=1}\Gamma^k_{ik}+\frac{1}{2}(\frac{\partial
ln\rho^n}{\partial x_i}+m\frac{\partial ln\rho^n}{\partial
x_i}-\frac{\partial ln\rho^n}{\partial
x_i})+\frac{n}{2}\frac{\partial ln\rho^{-m}}{\partial
x_i}\nonumber\\&=&\sum^m_{k=1}\Gamma^k_{ik}=0
\end{eqnarray}
Similar,by (3.9) and (3.12),at point $P$,
\begin{eqnarray}
\sum^{m+1}_{A=1}\tilde{\Gamma}^A_{\alpha
A}&=&\sum^{m+n}_{\beta=m+1}\Gamma^{\beta}_{\alpha\beta}=0
\end{eqnarray}
In what follows,compute the Ricci curvature of point $P$.Denote
$\tilde{R}_{AB},R_{ij},R_{\alpha\beta}$ the Ricci curvature of
$K,g,h$,
\begin{eqnarray}
\tilde{R}_{AB}=\sum^{m+n}_{C=1}\frac{\partial
\tilde{\Gamma}^C_{AB}}{\partial
x_C}-\sum^{m+n}_{C=1}\frac{\partial
\tilde{\Gamma}^C_{AC}}{\partial
x_B}+\sum^{m+n}_{C,D=1}\tilde{\Gamma}^D_{AB}\tilde{\Gamma}^C_{DC}-\sum^{m+n}_{C,D=1}
\tilde{\Gamma}^D_{AC}\tilde{\Gamma}^C_{BD}
\end{eqnarray}
First,compute $\tilde{R}_{ij}$,by (3.5),(3.7),at point $P$,
\begin{eqnarray}
\sum^{m+n}_{C=1}\frac{\partial \tilde{\Gamma}^C_{ij}}{\partial
x_C}&=&\sum^m_{k=1}\frac{\partial\tilde{\Gamma}^k_{ij}}{\partial
x_k}+\sum^{m+n}_{\alpha=m+1}\frac{\partial
\tilde{\Gamma}^{\alpha}_{ij}}{\partial
y_{\alpha}}\nonumber\end{eqnarray}\begin{eqnarray}&=&\sum^m_{k=1}\frac{\partial
\Gamma^k_{ij}}{\partial x_k}+\frac{\partial^2 ln\rho^n}{\partial
x_i\partial
x_j}-\frac{1}{2}g_{ij}\triangle_gln\rho^n\nonumber\\&&+\sum^{m+n}_{\alpha,\gamma=m+1}(-\frac{1}{2}h^{\alpha\gamma}g_{ij})
(\frac{\partial\rho^{m+n}}{\partial y_{\alpha}}\frac{\partial
ln\rho^n}{\partial y_{\gamma}}+\rho^{m+n}\frac{\partial^2
ln\rho^n}{\partial y_{\alpha}\partial
y_{\gamma}})\nonumber\\&=&\sum^k_{k=1}\frac{\partial\Gamma^k_{ij}}{\partial
x_k}+\frac{\partial^2 ln\rho^n}{\partial x_i\partial
x_j}-\frac{1}{2}g_{ij}\triangle_gln\rho^n-\frac{1}{2}g_{ij}\rho^{m+n}\triangle_hln\rho^n\nonumber\\&&
-\frac{1}{2}g_{ij}<\nabla_h\rho^{m+n}, \nabla_hln\rho^n>_h
\end{eqnarray}
Here,denote $\triangle_g,\triangle_h$ the Beltrami-Laplace operators
to $(M,g),(N,h)$ respectively.Use $\nabla_g,\nabla_h$ as the
covariant derivation to $M,N$,and $<\cdot,\cdot>_g,<\cdot,\cdot>_h$
as the metric $g,h$.And using (3.15),
\begin{equation}
\sum^{m+n}_{C=1}\frac{\partial\tilde{\Gamma}^C_{iC}}{\partial
x_j}=\sum^m_{k=1}\frac{\partial\Gamma^k_{ik}}{\partial x_j}
\end{equation}
By (3.15) and (3.16),at point $P$,
\begin{eqnarray}
\sum^{m+n}_{C,D=1}\tilde{\Gamma}^D_{ij}\tilde{\Gamma}^C_{DC}&=&\sum^m_{k=1}\sum^{m+n}_{C=1}\tilde{\Gamma}^k_{ij}
\tilde{\Gamma}^C_{kC}+\sum^{m+n}_{\alpha=m+1}\sum^{m+n}_{C=1}\tilde{\Gamma}^{\alpha}_{ij}\tilde{\Gamma}^C_{\alpha
C}=0\\\sum^{m+n}_{C,D=1}\tilde{\Gamma}^D_{iC}\tilde{\Gamma}^C_{jD}&=&\sum^m_{k,l=1}\tilde{\Gamma}^k_{il}\tilde{\Gamma}^l_{jk}
+\sum^m_{l=1}\sum^{m+n}_{\alpha=m+1}\tilde{\Gamma}^{\alpha}_{il}\tilde{\Gamma}^l_{j\alpha}+\sum^m_{k=1}\sum^{m+n}_{\alpha=m+1}
\tilde{\Gamma}^k_{i\alpha}\tilde{\Gamma}^{\alpha}_{jk}\nonumber\\&&+\sum^{m+n}_{\alpha,\beta=m+1}\tilde{\Gamma}^{\beta}_{i\alpha}\tilde{\Gamma}^{\alpha}_{j\beta}
\end{eqnarray}
Using (3.5),at point $P$ one has,
\begin{eqnarray}
\sum^m_{k,l=1}\tilde{\Gamma}^k_{il}\tilde{\Gamma}^l_{jk}&=&\sum^m_{k,l=1}\frac{1}{2}(\frac{\partial
ln\rho^n}{\partial x_i}\delta_{lk}+\frac{\partial ln\rho^n}{\partial
x_l}\delta_{ik}-\sum^m_{s=1}\frac{\partial ln\rho^n}{\partial
x_s}\delta_{sk}\delta_{il})\times
\nonumber\\&&\frac{1}{2}(\frac{\partial ln\rho^n}{\partial
x_j}\delta_{kl}+\frac{\partial ln\rho^n}{\partial
x_k}\delta_{jl}-\sum^m_{t=1}\frac{\partial ln\rho^n}{\partial
x_t}\delta_{tl}\delta_{jk})\nonumber\end{eqnarray}\begin{eqnarray}&=&\frac{1}{4}\sum^m_{k,l=1}[\frac{\partial
ln\rho^n}{\partial x_i}\frac{\partial ln\rho^n}{\partial
x_j}\delta_{kl}\delta_{kl}+\frac{\partial ln\rho^n}{\partial
x_l}\frac{\partial ln\rho^n}{\partial
x_j}\delta_{ik}\delta_{kl}-\frac{\partial ln\rho^n}{\partial
x_k}\frac{\partial ln\rho^n}{\partial
x_j}\delta_{il}\delta_{kl}\nonumber\\&&+\frac{\partial
ln\rho^n}{\partial x_i}\frac{\partial ln\rho^n}{\partial
x_k}\delta_{lk}\delta_{jl}+\frac{\partial ln\rho^n}{\partial
x_l}\frac{\partial ln\rho^n}{\partial
x_k}\delta_{ik}\delta_{jl}-\frac{\partial ln\rho^n}{\partial
x_k}\frac{\partial ln\rho^n}{\partial
x_k}\delta_{il}\delta_{jl}\nonumber\\&&-\frac{\partial
ln\rho^n}{\partial x_l}\frac{\partial ln\rho^n}{\partial
x_i}\delta_{lk}\delta_{jk}-\frac{\partial ln\rho^n}{\partial
x_l}\frac{\partial ln\rho^n}{\partial
x_l}\delta_{ik}\delta_{jk}+\frac{\partial ln\rho^n}{\partial
x_k}\frac{\partial ln\rho^n}{\partial
x_l}\delta_{il}\delta_{jk}]\nonumber\\&=&\frac{1}{4}[(m+2)\frac{\partial
ln\rho^n}{\partial x_i}\frac{\partial ln\rho^n}{\partial
x_j}-2|\nabla_gln\rho^n|^2_g\delta_{ij}]
\end{eqnarray}
Here,to a smooth enough function $u$,
$|\nabla_gu|^2_g=\sum^m_{i,j=1}g^{ij}\nabla_iu\nabla_ju$.Using (3.7)
and (3.9),at point $P$,
\begin{eqnarray}
\sum^m_{l=1}\sum^{m+n}_{\alpha=m+1}\tilde{\Gamma}^{\alpha}_{il}\tilde{\Gamma}^l_{j\alpha}&=&\sum^m_{l=1}
\sum^{m+n}_{\alpha,\beta=m+1}\rho^mh^{\alpha\beta}(-\frac{1}{2})\frac{\partial
\rho^n}{\partial
y_{\beta}}g_{il}\frac{1}{2}\delta_{jl}\frac{\partial
ln\rho^n}{\partial
y_{\alpha}}\nonumber\\&=&-\frac{1}{4}\rho^{m+n}\delta_{ij}|\nabla_hln\rho^n|^2_h
\end{eqnarray}
Here,$|\nabla_hu|^2_h=\sum^{m+n}_{\alpha,\beta=m+1}h^{\alpha\beta}\nabla_{\alpha}u\nabla_{\beta}u$,so,at
point $P$,
\begin{eqnarray}
\sum^m_{k=1}\sum^{m+n}_{\alpha=m+1}\tilde{\Gamma}^k_{i\alpha}\tilde{\Gamma}^{\alpha}_{jk}=-\frac{1}{4}\rho^{m+n}
\delta_{ij}|\nabla_h ln\rho^n|^2_h
\end{eqnarray}
By (3.11),at point $P$,
\begin{eqnarray}
\sum^{m+n}_{\alpha,\beta=m+1}\tilde{\Gamma}^{\beta}_{i\alpha}\tilde{\Gamma}^{\alpha}_{j\beta}&=&
\sum^{m+n}_{\alpha,\beta=m+1}\frac{1}{2}\delta_{\alpha\beta}\frac{\partial
ln\rho^{-m}}{\partial
x_i}\frac{1}{2}\delta_{\alpha\beta}\frac{\partial
ln\rho^{-m}}{\partial x_j}\nonumber\\&=&\frac{n}{4}\frac{\partial
ln\rho^{-m}}{\partial x_i}\frac{\partial ln\rho^{-m}}{\partial x_j}
\end{eqnarray}
Insert (3.22)-(3.25) to (3.21),at point $P$,
\begin{eqnarray}
\sum^{m+n}_{C,D=1}\tilde{\Gamma}^D_{iC}\tilde{\Gamma}^C_{jD}&=&\frac{1}{4}\{(m+2)\frac{\partial
ln\rho^n}{\partial x_i}\frac{\partial ln\rho^n}{\partial
x_j}-2|\nabla_g
ln\rho^n|^2_g\delta_{ij}\nonumber\\&&-2\rho^{m+n}\delta_{ij}|\nabla_hln\rho^n|^2_h+n\frac{\partial
ln\rho^{-m}}{\partial x_i}\frac{\partial ln\rho^{-m}}{\partial
x_j}\}\nonumber\\&=&\frac{1}{4}\{n(nm+2n+m^2)\frac{\partial
ln\rho}{\partial x_i}\frac{\partial ln\rho}{\partial
x_j}-2|\nabla_gln\rho^n|^2_g\delta_{ij}\nonumber\\&&-2\rho^{m+n}|\nabla_hln\rho^n|^2_h\delta_{ij}\}
\end{eqnarray}
By (3.17)-(3.20) and (3.26),at point $P$,
\begin{eqnarray}
\tilde{R}_{ij}&=&\sum^m_{k=1}\frac{\partial \Gamma^k_{ij}}{\partial
x_k}-\sum^m_{k=1}\frac{\partial \Gamma^k_{ik}}{\partial
x_j}+\frac{\partial^2ln\rho^n}{\partial x_i\partial
x_j}-\frac{1}{2}g_{ij}\triangle_gln\rho^n-\frac{1}{2}g_{ij}\rho^{m+n}\triangle_hln\rho^n\nonumber\\&&
-\frac{1}{2}g_{ij}<\nabla_h\rho^{m+n},\nabla_hln\rho^n>_h-\frac{1}{4}\{n(nm+2n+m^2)\frac{\partial
ln\rho}{\partial x_i}\frac{\partial ln\rho}{\partial
x_j}\nonumber\\&&-2|\nabla_g
ln\rho^n|^2_g\delta_{ij}-2\rho^{m+n}|\nabla_hln\rho^n|^2_h\delta_{ij}\}\nonumber\end{eqnarray}\begin{eqnarray}
&=&R_{ij}+\frac{\partial^2ln\rho^n}{\partial x_i\partial
x_j}-\frac{1}{2}g_{ij}\triangle_gln\rho^n-\frac{1}{2}g_{ij}\rho^{m+n}\triangle_hln\rho^n
\nonumber\\&&-\frac{m+n}{2}g_{ij}\rho^{m+n}n|\nabla_hln\rho|^2_h-\frac{1}{4}n(nm+2n+m^2)\frac{\partial
ln\rho}{\partial x_i}\frac{\partial ln\rho}{\partial
x_j}\nonumber\\&&+\frac{n^2}{2}|\nabla_gln\rho|^2_g\delta_{ij}+\frac{n^2}{2}|\nabla_hln\rho|^2_h\delta_{ij}\rho^{m+n}
\end{eqnarray}
Using (3.2) and (3.27),at point $P$ one has,
\begin{eqnarray}
\sum^m_{i,j=1}K^{ij}\tilde{R}_{ij}&=&\rho^{-n}R_g+\rho^{-n}\triangle_gln\rho^n-\frac{m}{2}\rho^{-n}\triangle_gln\rho^n
-\frac{m}{2}\rho^m\triangle_hln\rho^n\nonumber\\&&-\frac{mn}{2}m\rho^m|\nabla_hln\rho|^2_h-\frac{n}{4}\rho^{-n}(nm+2n+m^2)|\nabla_gln\rho|^2_g
\nonumber\\&&+\frac{mn^2}{2}|\nabla_gln\rho|^2_g\rho^{-n}\nonumber\\&=&
\rho^{-n}R_g+\frac{2-m}{2}\rho^{-n}\triangle_gln\rho^n-\frac{m}{2}\rho^m\triangle_hln\rho^n\nonumber\\&&
-\frac{n}{4}(-nm+2n+m^2)\rho^{-n}|\nabla_gln\rho|^2_g\nonumber\\&&-\frac{m^2n}{2}\rho^m|\nabla_hln\rho|^2_h
\end{eqnarray}
Here,$R_g,R_h$ denote the scalar curvature of $(M,g),(N,h)$.Since
the position of $g$ and $h$ are symmetry,interchange $g$ and
$h$,$\rho$ and $\rho^{-1}$,$m$ and $n$ of (3.28) at point $P$,one
has,
\begin{eqnarray}
\sum^{m+n}_{\alpha,\beta=m+1}K^{\alpha\beta}\tilde{R}_{\alpha\beta}&=&\rho^mR_h+\frac{2-n}{2}(\triangle_hln\rho^{-m})\rho^m-\frac{n}{2}\rho^{-n}\triangle_gln\rho^{-m}
\nonumber\\&&-\frac{m}{4}(-nm+2m+n^2)\rho^m|\nabla_hln\rho|^2_h\nonumber\\&&-\frac{n^2m}{2}\rho^{-n}|\nabla_gln\rho|^2_g
\end{eqnarray}
Take $\tilde{R}$ as the scalar curvature of $(\Sigma,K)$,then at
point $P$,by (3.28),(3.29) and
\begin{equation}
\tilde{R}=\sum^m_{i,j=1}K^{ij}\tilde{R}_{ij}+\sum^{m+n}_{\alpha,\beta=m+1}K^{\alpha\beta}\tilde{R}_{\alpha\beta}
\end{equation}
one has (1.1),and let $\rho=e^{2u}$, where $u$ is a smooth enough
function on $\Sigma$.Then one has,
\begin{eqnarray}
\tilde{R}&=&e^{-2nu}R_g+e^{2mu}R_h+2ne^{-2nu}\triangle_gu-2me^{2mu}\triangle_hu\nonumber\\&&
-n(nm+2n+m^2)e^{-2nu}|\nabla_gu|^2_g\nonumber\\&&
-m(nm+2m+n^2)e^{2mu}|\nabla_hu|^2_h
\end{eqnarray}
Since,the left and right hand sides of the above equation both are
not depended on the choice of local coordinates on $\Sigma$,then it
is a equation on manifold.If $n=1$,simplify it,then one obtains
\begin{eqnarray}
\tilde{R}&=&e^{-2u}R_g+e^{2mu}R_h+2e^{-2u}\triangle_gu-2me^{2mu}\triangle_hu\nonumber\\&&
-(m^2+m+2)e^{-2u}|\nabla_gu|^2_g-m(3m+1)e^{2mu}|\nabla_hu|^2_h
\end{eqnarray}
Since $n=1,N$ becomes a curve.Suppose $t$ be $N$'s arch length
parameter,then write $u$ as $u(x,t)$.So with this parameter,one
has,
\begin{equation}
h(\frac{\partial}{\partial t},\frac{\partial}{\partial
t})=1;\triangle_hu=\frac{\partial^2u}{\partial
t^2};\nabla_hu=\frac{\partial u}{\partial t};R_h=0
\end{equation}
By (3.32)(3.33),one obtains (1.2).Here,we take the conservation that
$\triangle=\triangle_g,\nabla u=\nabla_gu,|\nabla
u|^2=|\nabla_gu|^2_g,<\cdot,\cdot>=<\cdot,\cdot>_g$.And
 the following paragraph holds the conservation.
 \section{Prior Estimate:Energy Estimate}
\quad\ The idea of the energy estimate used here is inspired by
[Smith$\&$Tataru].Let
\begin{eqnarray}
F(u,v)&=&\frac{1}{2}[-e^{-2mv}\tilde{R}+e^{-2(m+1)v}R_g\nonumber\\&&-(m^2+m+2)e^{-2(m+1)v}<\nabla
v,\nabla u>-m(3m+1)v_tu_t]
\end{eqnarray}
In this section,consider the prior estimate of the following
equation,
\begin{eqnarray}
mu_{tt}-e^{-2(m+1)v}\triangle u=F(u,v)
\end{eqnarray}
when $u=v$,(4.2) becomes to (1.2).Here,suppose $(M,g)$ be a closed
Riemainnian manifold of m-dimension.$I=[0,T]$ is a interval of
length T.And $u,v\in C^{s+1}(M\times I),s>\frac{m}{2}+1$.First
estimate
$\|me^{2(m+1)v}\nabla^{l-1}u_{tt}-\triangle\nabla^{l-1}u\|_{L^2(M)}$,here
$1\leq l\leq s$.Then,by(4.2)
\begin{eqnarray}
&&me^{2(m+1)v}\nabla^{l-1}u_{tt}-\triangle\nabla^{l-1}u\nonumber\\&=&[e^{2(m+1)v}\nabla^{l-1}F(u,v)]+[e^{2(m+1)v}\nabla^{l-1}(e^{-2(m+1)v}\triangle
u)-\nabla^{l-1}\triangle u]\nonumber\\
&&+[\nabla^{l-1}\triangle u-\triangle\nabla^{l-1}
u]\nonumber\\&=&I+II+III
\end{eqnarray}
In what follows,estimate $I,II,III$ respectively.
\begin{eqnarray}
\|I\|_{L^2(M)}&\leq&
Ce^{C\|v\|_{C^0(M)}}\|\nabla^{l-1}[-e^{-2mv}\tilde{R}+e^{-2(m+1)v}R_g\nonumber\\&&-(m^2+m+2)e^{-2(m+1)v}<\nabla
v,\nabla u>-m(3m+1)v_tu_t]\|_{L^2(M)}\nonumber\\&&
\end{eqnarray}
By Corollary 2.4 and Sobolev embedding theorem,
$H^{s-1}(M)\hookrightarrow C^0(M)$ one has,
\begin{eqnarray}
&&\|\nabla^{l-1}(e^{-2mv}\tilde{R})\|_{L^2(M)}\nonumber\\&\leq&\|\nabla^{l-1}[(e^{-2mv}-1)\tilde{R}]\|_{L^2(M)}
+\|\nabla^{l-1}\tilde{R}\|_{L^2(M)} \nonumber\\&\leq
&C(\|e^{-2mv}-1\|_{C^0(M)}\|\tilde{R}\|_{H^{l-1}(M)}+\|e^{-2mv}-1\|_{H^{l-1}(M)}\|\tilde{R}\|_{C^{0}(M)})
\nonumber\\&&+\|\nabla^{l-1}\tilde{R}\|_{L^2(M)}\nonumber\\&\leq&C\|e^{-2mv}-1\|_{H^{s-1}(M)}\|\tilde{R}\|_{H^{s-1}(M)}
+\|\nabla^{l-1}\tilde{R}\|_{L^2(M)}\nonumber\\&\leq&C\|e^{-2mw}-1\|_{C^{s-1}[-\|v\|_{C^0(M)},\|v\|_{C^0(M)}]}(1+\|v\|^{s-1}_{C^0(M)})
\|v\|_{H^{s-1}(M)}\|\tilde{R}\|_{H^{s-1}(M)}\nonumber\\&&+\|\nabla^{l-1}\tilde{R}\|_{L^2(M)}\nonumber\\
&\leq&Ce^{C\|v\|_{C^0(M)}}(1+\|v\|_{H^{s-1}(M)})\|\tilde{R}\|_{H^{s-1}(M)}\nonumber\\
&\leq&Ce^{C\|v\|_{H^{s-1}(M)}}\|\tilde{R}\|_{H^{s-1}(M)}
\end{eqnarray}
Similar,
\begin{eqnarray}
\|\nabla^{l-1}(e^{-2(m+1)v}R_g)\|_{L^2(M)}&\leq&Ce^{C\|v\|_{H^{s-1}(M)}}\|R_g\|_{H^{s-1}(M)}
\end{eqnarray}
Again,by Corollary 2.4 and sobolev embedding theorem,
$H^{s-1}(M)\hookrightarrow C^0(M)$, $H^{s}(M)\hookrightarrow C^1(M)$
one has,
\begin{eqnarray}
&&\|\nabla^{l-1}(e^{-2(m+1)v}<\nabla v,\nabla
u>)\|_{L^2(M)}\nonumber\\&=&C\|\nabla^{l-1}<\nabla
(e^{-2(m+1)v}-1),\nabla u>\|_{L^2(M)}\nonumber\\&\leq
&C(\|e^{-2(m+1)v}-1\|_{C^1(M)}\|u\|_{H^{l}(M)}
+\|e^{-2(m+1)v}-1\|_{H^{l}(M)}\|u\|_{C^1(M)})\nonumber\\&\leq&
C\|e^{-2(m+1)v}-1\|_{H^s(M)}\|u\|_{H^s(M)}\nonumber\\&\leq
&Ce^{\|v\|_{H^s(M)}}\|u\|_{H^s(M)}
\end{eqnarray}
the last step is similar to (4.5).And
\begin{eqnarray}
\|\nabla^{l-1}(v_tu_t)\|_{L^2(M)}&\leq&C(\|v_t\|_{H^{l-1}(M)}\|u_t\|_{C^0(M)}+\|v_t\|_{C^0(M)}\|u_t\|_{H^{l-1}(M)})\nonumber\\
&\leq&C\|v_t\|_{H^{s-1}(M)}\|u_t\|_{H^{s-1}(M)}
\end{eqnarray}
Conservation that,
\begin{eqnarray}
A_s=\|\tilde{R}\|_{H^{s-1}(M)}+\|R_g\|_{H^{s-1}(M)}
\end{eqnarray}
Then,by (4.4)-(4.9) and $H^{s-1}(M)\hookrightarrow C^0(M)$, one
has,
\begin{eqnarray}
\|I\|_{L^2(M)}\leq
Ce^{C(\|v\|_{H^s(M)}+\|v_t\|_{H^{s-1}(M)})}(A_s+\|u\|_{H^s(M)}+\|u_t\|_{H^{s-1}(M)})
\end{eqnarray}
For $II$,
\begin{eqnarray}
\|II\|_{L^2(M)}&=&\|e^{2(m+1)v}\nabla^{l-1}(e^{-2(m+1)v}\triangle
u)-\nabla^{l-1}\triangle
u\|_{L^2(M)}\nonumber\\&=&\|e^{2(m+1)v}\sum_{i+j=l-1,i\geq
1}\nabla^i(e^{-2(m+1)v})\nabla^j\triangle u\|_{L^2(M)}\nonumber\\
&\leq&Ce^{C\|v\|_{C^0(M)}}\sum_{i+j=l-1,i\geq
1}\|\nabla^i(e^{-2(m+1)v}-1)\nabla^j\triangle
u\|_{L^2(M)}\nonumber\\
\end{eqnarray}
For fixed $i,j$, let,
\begin{eqnarray}
p=\frac{l-1}{i-1},q=\frac{l-1}{j+1},
\frac{1}{2p}+\frac{1}{2q}=\frac{1}{2},
\end{eqnarray}
Then,
\begin{eqnarray}
\frac{1}{2p}-\frac{i-1}{m}=\frac{i-1}{l-1}(\frac{1}{2}-\frac{l-1}{m}),
\frac{1}{2q}-\frac{j+1}{m}=\frac{j+1}{l-1}(\frac{1}{2}-\frac{l-1}{m})
\end{eqnarray}
Using H\"{o}lder inequality,Corollary 2.4,Remark 2.6 and Sobolev embedding
theorem,one obtains,
\begin{eqnarray}
&&\|\nabla^i(e^{-2(m+1)v}-1)\nabla^j\triangle
u\|_{L^2(M)}\nonumber\\&\leq&\|\nabla^i(e^{-2(m+1)v}-1)\|_{L^{2p}(M)}\|\nabla^j\triangle
u\|_{L^{2q}(M)}\nonumber\\&\leq&C\|\nabla^i(e^{-2(m+1)v}-1)\|_{L^{2p}(M)}\|\nabla^{j+2}
u\|_{L^{2q}(M)}\nonumber\\&\leq&C\|e^{-2(m+1)v}-1\|^{\frac{i-1}{l-1}}_{H^{l}(M)}
\|e^{-2(m+1)v}-1\|^{\frac{l-i}{l-1}}_{C^{1}(M)}\|u\|^{\frac{j+1}{l-1}}_{H^{l}(M)}\|u\|^{\frac{l-j-2}{l-1}}_{C^1(M)}\nonumber\\
&\leq&C\|e^{-2(m+1)v}-1\|^{\frac{i-1}{l-1}}_{H^{s}(M)}
\|e^{-2(m+1)v}-1\|^{\frac{l-i}{l-1}}_{H^s(M)}\|u\|^{\frac{j+1}{l-1}}_{H^{s}(M)}\|u\|^{\frac{l-j-2}{l-1}}_{H^s(M)}\nonumber\\
&\leq&C\|e^{-2(m+1)v}-1\|_{H^{s}(M)}\|u\|_{H^{s}(M)}\nonumber\\&\leq&Ce^{C\|v\|_{C^0(M)}}\|v\|_{H^s(M)}\|u\|_{H^s(M)}
\end{eqnarray}
Then,
\begin{eqnarray}
\|II\|_{L^2(M)}&\leq&Ce^{C\|v\|_{C^0(M)}}\|v\|_{H^s(M)}\|u\|_{H^s(M)}
\end{eqnarray}
For $III$,denote $B_{i_1\cdots i_k}=\nabla_{i_k}\cdots \nabla_{i_1}
u$, by the Ricci identity,
\begin{eqnarray}
(\nabla_p\nabla_q-\nabla_q\nabla_p)B_{i_1\cdots
i_k}=-\sum^k_{t=1}\sum^m_{n=1}R^n_{pqi_t}B_{i_1\cdots
i_{t-1}ni_{t+1}\cdots i_k}
\end{eqnarray}
here,all index take values in $\{1,\cdots,m\}$,then direct
computation gives,
\begin{eqnarray}
&&\nabla_j\triangle B_{i_1\cdots
i_k}\nonumber\\&=&\triangle\nabla_jB_{i_1\cdots
i_k}-\sum^m_{p,q,n=1}\sum^k_{t=1}g^{pq}\nabla_pR^n_{jqi_t}B_{i_1\cdots
i_{t-1}ni_{t+1}\cdots
i_k}\nonumber\\&&-\sum^m_{p,q,n=1}\sum^k_{t=1}g^{pq}R^n_{jqi_t}\nabla_pB_{i_1\cdots
i_{t-1}ni_{t+1}\cdots
i_k}-\sum^m_{p,q,n=1}\sum^k_{t=1}g^{pq}R^n_{jpi_t}\nabla_qB_{i_1\cdots
i_{t-1}ni_{t+1}
i_k}\nonumber\\&&-\sum^m_{p,q,n=1}g^{pq}R^n_{jpq}\nabla_nB_{i_1\cdots
i_k}
\end{eqnarray}
From (4.17),one can obtain,if $0\leq \tilde{k}$,
\begin{eqnarray}
\|\nabla^{\tilde{k}}(\nabla\triangle B-\triangle\nabla
B)\|_{L^2(M)}&\leq& C\|u\|_{H^{k+\tilde{k}+1}(M)}
\end{eqnarray}
here, $C$ is a positive constant depended only on manifold
$M$,then iterate (4.18),one has,
\begin{eqnarray}
\|\nabla\triangle\nabla^{l-2}u-\triangle\nabla^{l-1}u\|_{L^2(M)}&\leq&
C\|u\|_{H^{l-1}(M)}\nonumber\\
\|\nabla^2\triangle\nabla^{l-3}u-\nabla\triangle\nabla^{l-2}u\|_{L^2(M)}&\leq&
C\|u\|_{H^{l-1}(M)}\nonumber\\\cdots\cdots&&\nonumber\\
\|\nabla^{l-1}\triangle u-\nabla^{l-2}\triangle u\|_{L^2(M)}&\leq&
C\|u\|_{H^{l-1}(M)}
\end{eqnarray}
adding the left and right hand sides of (4.19),then one obtains
\begin{eqnarray}
\|III\|_{L^2(M)}\leq C\|u\|_{H^s(M)}
\end{eqnarray}
By (4.3)(4.9)(4.10)(4.15)(4.20) and Sobolev embedding theorem one
obtains,
\begin{eqnarray}
&&\|me^{2(m+1)v}\nabla^{l-1}u_{tt}-\triangle\nabla^{l-1}
u\|_{L^2(M)}\nonumber\\&\leq&
C_1e^{C_2(\|v\|_{H^s(M)}+\|v_t\|_{H_{s-1}(M)})}(A_s+\|u\|_{H^s(M)}+\|u_t\|_{H^{s-1}(M)})
\end{eqnarray}
here,$C_1,C_2$ are positive constants depended on $s$,manifold
$M$,not depended on $\tilde{R},u$.

Then we do energy estimate.One can introduce the following energy
functions:
\begin{eqnarray}
E^{(l)}_1(t)&=&\frac{1}{2}\int_M[|\nabla^lu|^2+me^{2(m+1)v}|\nabla^{l-1}u_t|^2] dV \\
E^{(l)}_2(t)&=&\frac{1}{2}\int_M[e^{-2(m+1)v}|\nabla^lu|^2+m|\nabla^{l-1}u_t|^2]dV
\end{eqnarray}
and let,
\begin{eqnarray}
E_s(t)=\sum^s_{l=1}(E^{(l)}_1(t)+E^{(l)}_2(t))+\frac{1}{2}\int_M|u|^2dV
\end{eqnarray}
Obvious,
\begin{eqnarray}
\|u\|_{H^s(M)}+\|u_t\|_{H^{s-1}(M)}\leq C(E_s(t))^{\frac{1}{2}}
\end{eqnarray}
To (4.22),take derivation respect to $t$,by (2.2) one has,
\begin{eqnarray}
\frac{d}{dt}E^{(l)}_1(t)&=&\int_M[<\nabla^lu_t,\nabla^l
u>+m(m+1)e^{2(m+1)v}v_t|\nabla^{l-1}u_t|^2\nonumber\\&&+me^{2(m+1)v}<\nabla^{l-1}u_t,\nabla^{l-1}u_{tt}>]dV\nonumber\\
&\leq&C\|v_t\|_{C^0(M)}E^{(l)}_1(t)+\int_M[<\nabla^lu_t,\nabla^lu>\nonumber\\&&+<\nabla^{l-1}u_t,me^{2(m+1)v}\nabla^{l-1}u_{tt}>]dV
\end{eqnarray}
Since $M$ is closed,$\int_M
div<\nabla^{l-1}u_t,\nabla\nabla^{l-1}u>dV=0$.This implies,
\begin{eqnarray}
\int_M<\nabla^lu_t,\nabla^lu>dV=-\int_M<\nabla^{l-1}u_t,\triangle\nabla^{l-1}u>dV
\end{eqnarray}
then,by (4.26),and H\"{o}lder inequality,
\begin{eqnarray}
&&\frac{d}{dt}E^{(l)}_1(t)\nonumber\\&\leq&C\|v_t\|_{C^0(M)}E^{(l)}_1(t)+\int_M<\nabla^{l-1}u_t,me^{2(m+1)v}\nabla^{l-1}u_{tt}
-\triangle\nabla^{l-1}u>dV\nonumber\\&\leq&C\|v_t\|_{C^0(M)}E^{(l)}_1(t)+\|\nabla^{l-1}u_t\|_{L^2(M)}\|me^{2(m+1)v}\nabla^{l-1}u_{tt}-\triangle\nabla^{l-1}
u\|_{L^2(M)}\nonumber\\&\leq&C\|v_t\|_{C^0(M)}E^{(l)}_1(t)\nonumber\\
&&+Ce^{C\|v\|_{C^0(M)}}\|me^{2(m+1)v}\nabla^{l-1}u_{tt}-\triangle\nabla^{l-1}u\|_{L^2(M)}(E^{(l)}_1(t))^{\frac{1}{2}}
\end{eqnarray}
Similar,
\begin{eqnarray}
&&\frac{d}{dt}E^{(l)}_2(t)\nonumber\\&=&\int_M[-e^{-2(m+1)v}(m+1)v_t|\nabla^lu|^2+<e^{-2(m+1)v}\nabla^lu,\nabla^lu_t>
\nonumber\\&&+<\nabla^{l-1}u_t,m\nabla^{l-1}u_{tt}>]dV\nonumber\\&\leq&
C\|v_t\|_{C^0(M)}E^{(l)}_2(t)\nonumber\\&&+\int_M[<\nabla^{l-1}u_t,m\nabla^{l-1}u_{tt}>
-<\nabla^{l-1}u_t,e^{-2(m+1)v}\triangle\nabla^{l-1}u>]dV\nonumber\\&&+2(m+1)\int_Me^{-2(m+1)v}<\sum^m_{i,j=1}g^{ij}
\nabla_iv\nabla_j\nabla^{l-1}u,\nabla^{l-1}u_t>dV \nonumber\\ &\leq&
C\|v_t\|_{C^0(M)}E^{(l)}_2(t)+\|\nabla^{l-1}u_t\|_{L^2(M)}\|m\nabla^{l-1}u_{tt}-e^{-2(m+1)v}\triangle\nabla^{l-1}u\|_{L^2(M)}
\nonumber\\&&+Ce^{C\|v\|_{C^0(M)}}\|\sum^m_{i,j=1}g^{ij}\nabla_iv\nabla_j\nabla^{l-1}u\|_{L^2(M)}\|\nabla^{l-1}u_t\|_{L^2(M)}
\nonumber\\&\leq&C\|v_t\|_{C^0(M)}E^{(l)}_2(t)\nonumber\\&&+Ce^{C\|v\|_{C^0(M)}}\|me^{2(m+1)v}\nabla^{l-1}u_{tt}-\triangle\nabla^{l-1}u\|_{L^2(M)}(E^{(l)}_2(t))^{\frac{1}{2}}
\nonumber\\&&+Ce^{C\|v\|_{C^0(M)}}\|\nabla
v\|_{C^0(M)}\|\nabla^lu\|_{L^2(M)}\|\nabla^{l-1}u_t\|_{L^2(M)}\nonumber\\
&\leq& Ce^{C\|v\|_{C^0(M)}}(\|v_t\|_{C^0(M)}+\|\nabla
v\|_{C^0(M)})E^{(l)}_2(t)\nonumber\\&&
+Ce^{C\|v\|_{C^0(M)}}\|me^{2(m+1)v}\nabla^{l-1}u_{tt}-\triangle\nabla^{l-1}u\|_{L^2(M)}(E^{(l)}_2(t))^{\frac{1}{2}}
\end{eqnarray}
Using (4.28)(4.29),one has
\begin{eqnarray}
\frac{d}{dt}E_s(t)&\leq&
Ce^{C\|v\|_{C^0(M)}}(\|v_t\|_{C^0(M)}+\|\nabla
v\|_{C^0(M)})E_s(t)\nonumber\\&&+Ce^{C\|v\|_{C^0(M)}}\sum^s_{l=1}\|me^{2(m+1)v}\nabla^{l-1}u_{tt}-\triangle\nabla^{l-1}u\|_{L^2(M)}(E_s(t))^{\frac{1}{2}}
\nonumber\\&&+\int_Muu_tdV
\end{eqnarray}
By H\"{o}lder inequality,
\begin{eqnarray}
\int_Muu_tdV\leq CE_s(t)
\end{eqnarray}
then by (4.21)(4.25)(4.30)(4.31) and Sobolev embedding theorem,one
obtains
\begin{eqnarray}
&&\frac{d}{dt}E_s(t)\nonumber\\&\leq&Ce^{C\|v\|_{C^0(M)}}(1+\|v_t\|_{C^0(M)}+\|\nabla
v\|_{C^0(M)})E_s(t)\nonumber\\&&+Ce^{C(\|v\|_{H^s(M)}+\|v_t\|_{H^{s-1}(M)})}(A_s+\|u\|_{H^s(M)}+\|u_t\|_{H^{s-1}(M)})
(E_s(t))^{\frac{1}{2}}\nonumber\\
&\leq&Ce^{C(\|v\|_{H^s(M)}+\|v_t\|_{H^{s-1}(M)})}(A_s+\|u\|_{H^s(M)}+\|u_t\|_{H^{s-1}(M)}+(E_s(t))^{\frac{1}{2}})(E_s(t))^{\frac{1}{2}}
\nonumber\\&\leq&Ce^{C(\|v\|_{H^s(M)}+\|v_t\|_{H^{s-1}(M)})}(A_s+(E_s(t))^{\frac{1}{2}})(E_s(t))^{\frac{1}{2}}
\end{eqnarray}
Then obvious,
\begin{eqnarray}
\frac{d}{dt}(E_s(t))^{\frac{1}{2}}&\leq&C_1e^{C_2(\|v\|_{H^s(M)}+\|v_t\|_{H^{s-1}(M)})}(A_s+(E_s(t))^{\frac{1}{2}})
\end{eqnarray}
here,$C_1,C_2$ related to manifold $M$ and $s$.Then using (4.33),one has the following Proposition.\\

\noindent\textbf{Proposition 4.1} \quad\ Suppose $(M,g)$ be
m-dimensional closed Reimainnian manifold,and $I=[0,T]$ be a
interval of length $T$.Let $u,v\in C^{s+1}(M\times
I)$,$s>\frac{m}{2}+1$,and $u,v$ satisfy
(4.2) then:\\
(i)If exist a positive constant
$D>2\sqrt{2}[1+(E_s(0))^{\frac{1}{2}}]$.There is $t_0>0$ which is
related to manifold
$M$,$s,\|\tilde{R}\|_{C^0(I,H^{s-1}(M))},D,E_s(0)$, such that
\begin{eqnarray}
\|v\|_{H^s(M)}+\|v_t\|_{H^{s-1}(M)}&\leq& D
\end{eqnarray}
held in $0\leq t\leq t_0$ implies
\begin{eqnarray}
\|u\|_{H^s(M)}+\|u_t\|_{H^{s-1}(M)}&\leq& D
\end{eqnarray}
held in $0\leq t\leq t_0$\\
(ii)If $u=v$ which means $u$ satisfy (1.2),then exist
$t_0,\tilde{C}>0$ such that
\begin{eqnarray}
\|u\|_{H^s(M)}+\|u_t\|_{H^{s-1}(M)}&\leq& \tilde{C}
\end{eqnarray}
holds in $0\leq t\leq t_0$.Here,$t_0,\tilde{C}$ relates to manifold
$M$,$s,\|\tilde{R}\|_{C^0(I,H^{s-1}(M))},E_s(0)$, and $\tilde{C}$
also relates to $t_0$.\\
\textbf{Proof} \quad\ (i) Using (4.9) and (4.33),
\begin{eqnarray}
\frac{d}{dt}(E_s(t))^{\frac{1}{2}}&\leq&C_3e^{C_2(\|v\|_{H^s(M)}+\|v_t\|_{H^{s-1}(M)})}(1+(E_s(t))^{\frac{1}{2}})
\end{eqnarray}
here, $C_3$ related to manifold $M$,$s$,and extra
$\|\tilde{R}\|_{C^0(I,H^{s-1}(M))}$.Take $t_0$, such that
\begin{eqnarray}
2\sqrt{2}[1+(E_s(0))^{\frac{1}{2}}]e^{C_3e^{C_2D}t_0}&\leq& D
\end{eqnarray}
Then in $0\leq t\leq t_0$,
\begin{eqnarray}
\frac{d}{dt}(E_s(t))^{\frac{1}{2}}&\leq&C_3e^{C_2D}(1+(E_s(t))^{\frac{1}{2}})
\end{eqnarray}
Integrate (4.39) one has,
\begin{eqnarray}
(E_s(t))^{\frac{1}{2}}\leq
[1+(E_s(0))^{\frac{1}{2}}]e^{C_3e^{C_2D}t}\leq \frac{D}{2\sqrt{2}}
\end{eqnarray}
This implies the result by (4.25).\\
(ii)Let $v=u$ in (4.37),then
\begin{eqnarray}
\frac{d}{dt}(E_s(t))^{\frac{1}{2}}&\leq&C_3e^{C_2(E_s(t))^{\frac{1}{2}}}
\end{eqnarray}
Integrate it,
\begin{eqnarray}
e^{-C_2(E_s(t))^{\frac{1}{2}}}&\geq&e^{-C_2(E_s(0))^{\frac{1}{2}}}-C_3C_2t
\end{eqnarray}
Hence,when $t_0<\frac{e^{-C_2(E_s(0))^{\frac{1}{2}}}}{C_3C_2}$,one
has,
\begin{eqnarray}
(E_s(t))^{\frac{1}{2}}&\leq&\frac{1}{C_2}
ln\frac{1}{e^{-C_2(E_s(0))^{\frac{1}{2}}}-C_3C_2t_0}
\end{eqnarray}
This gives what we want.\\

\noindent\textbf{Proposition 4.2}\quad\ Suppose $(M,g)$ be
m-dimensional closed Reimainnian manifold,whose scalar curvature
is zero,and $I=[0,T]$ be a interval of length $T$.Let $u,v\in
C^{s+1}(M\times I)$, $s>\frac{m}{2}+1$. $u,v$ satisfy (4.2),and
$u(x,0)=u_t(x,0)=0$ to any $x\in M$.Then to any given positive
constant $D$,exist constant $\varepsilon >0$ which depends on
manifold $M,s,D$, such that
\begin{eqnarray}
\|v\|_{H^s(M)}+\|v_t\|_{H^{s-1}(M)}&\leq& D
\end{eqnarray}
held in $I$ implies
\begin{eqnarray}
\|u\|_{H^s(M)}+\|u_t\|_{H^{s-1}(M)}&\leq& D
\end{eqnarray}
held in $I$ when $\|\tilde{R}\|_{H^{s-1}(M)}\leq\varepsilon$.\\
\textbf{Proof}\quad\ Take $\varepsilon$ satisfying
\begin{eqnarray}
\varepsilon e^{C_1e^{C_2D}T}&\leq&\frac{D}{2\sqrt{2}}
\end{eqnarray}
here,$C_1,C_2$ are constant of (4.33).Since
$R_g=0$,$A_s=\|\tilde{R}\|_{H^{s-1}(M)}$, so by (4.33) and (4.44),
\begin{eqnarray}
\frac{d}{dt}(E_s(t))^{\frac{1}{2}}&\leq&C_1e^{C_2D}(\varepsilon+(E_s(t))^{\frac{1}{2}})
\end{eqnarray}
Integrate (4.47)
\begin{eqnarray}
(E_s(t))^{\frac{1}{2}}&\leq&\varepsilon e^{C_1e^{C_2D}T}\leq
\frac{D}{2\sqrt{2}}
\end{eqnarray}
and using (4.25),then one gets the result.
\section{Linear Equation}
\quad\ In this section we consider the linear equation (1.3) on
the condition:
\begin{eqnarray}
 u(x,0)=\varphi(x),u_t(x,0)=\psi(x)
\end{eqnarray}
here, $(M,g)$ the m-dimensional closed Reimainnian manifold,
$I=[0,T]$ interval of length $T$,and $a,\alpha,\beta,\gamma,f\in
C^{\infty}(M\times I),\varphi,\psi\in C^{\infty}(M)$, and exist
positive constant $L$,such that $\frac{1}{L}\leq\alpha\leq L$.In
what follows, we discuss the existence and property of the
solution.If $M$ is a domain of $\mathbb{R}^n$, then the following
results are classical in linear hyperbolic
equation(c.f.[H\"{o}mander] and [Xin]).We follow their method.
First,give some prior estimates. Multiple  $u_t$ in
 (1.3),and integrate on $M$.Since $M$ is closed, one has,
\begin{eqnarray}
&&\frac{1}{2}\frac{d}{dt}(\|u_t\|^2_{L^2(M)}+\int_M\alpha|\nabla
u|^2dV)\nonumber\\&\leq&\|f\|_{L^2(M)}\|u_t\|_{L^2(M)}-\int_Mu_t<\nabla(\alpha-\beta),\nabla
u>dV+C\int_M|\nabla u|^2dV\nonumber\\&&+
C\|u\|_{L^2(M)}\|u_t\|_{L^2(M)}+C\|u_t\|^2_{L^2(M)}
\end{eqnarray}
Using,
\begin{eqnarray}
\|u\|^2_{L^2(M)}&\leq&\int^t_0\|u_t\|^2_{L^2(M)}d\tau+\|\varphi\|^2_{L^2(M)}
\end{eqnarray}
one has,
\begin{eqnarray}
&&\|u_t\|^2_{L^2(M)}+\|\nabla
u\|^2_{L^2(M)}\nonumber\\&\leq&C[\|\varphi\|^2_{H^1(M)}+\|\psi\|^2_{L^2(M)}+\int^t_0(\|u_t\|^2_{L^2(M)}+\|\nabla
u\|^2_{L^2(M)}+\|f\|^2_{L^2(M)})d\tau]\nonumber\\
\end{eqnarray}
Integrate it,one has
\begin{eqnarray}
\|u_t\|^2_{L^2(M)}+\|\nabla
u\|^2_{L^2(M)}&\leq&C[\|\varphi\|^2_{H^1(M)}+\|\psi\|^2_{L^2(M)}+\int^t_0\|f\|^2_{L^2(M)}d\tau]\nonumber\\
\end{eqnarray}
Multiple $-\triangle u_t$ in (1.3),integrate on $M$,and using
(5.3)(5.5) $W^{k,p}$ estimate of elliptic equation,one has,
\begin{eqnarray}
&&\|\nabla u_t\|^2_{L^2(M)}+\frac{1}{L}\|\triangle
u\|^2_{L^2(M)}+2\int_M(\triangle u)
fdV\nonumber\\&\leq&C[\|\varphi\|^2_{H^2(M)}+\|\psi\|^2_{H^1(M)}+\|f(0)\|^2_{L^2(M)}+\int^t_0(\|f_t\|^2_{L^2(M)}+\|f\|^2_{L^2(M)}
\nonumber\\&&+\|\triangle u\|^2_{L^2(M)}+\|\nabla
u_t\|^2_{L^2(M)})d\tau]
\end{eqnarray}
Using the inequality $2\int_M(\triangle u) f\leq
\frac{1}{2L}\|\triangle u\|^2_{L^2(M)}+C\|f\|^2_{L^2(M)} $
and,
\begin{eqnarray}
\|f(t)\|^2_{L^2(M)}&\leq&\int^t_0\|f_t(\tau)\|^2_{L^2(M)}d\tau+\|f(0)\|^2_{L^2(M)}
\end{eqnarray}
One has,
\begin{eqnarray}
&&\|\nabla u_t\|^2_{L^2(M)}+\|\triangle u\|^2_{L^2(M)}
\nonumber\\&\leq&C[\|\varphi\|^2_{H^2(M)}+\|\psi\|^2_{H^1(M)}+\|f(0)\|^2_{L^2(M)}+\int^t_0(\|f_t\|^2_{L^2(M)}
+\|\triangle u\|^2_{L^2(M)}\nonumber\\&&+\|\nabla
u_t\|^2_{L^2(M)})d\tau]
\end{eqnarray}
Integrate it,
\begin{eqnarray}
\|\nabla u_t\|^2_{L^2(M)}+\|\triangle
u\|^2_{L^2(M)}&\leq&C[\|\varphi\|^2_{H^2(M)}+\|\psi\|^2_{H^1(M)}+\|f(0)\|^2_{L^2(M)}\nonumber\\&&+\int^t_0\|f_t\|^2_{L^2(M)}
d\tau]
\end{eqnarray}
Take derivation to $t$ in (1.3),and multiple $u_{tt}$ in
it,integrate on $M$,using (5.3)(5.5)(5.7)  one has,
\begin{eqnarray}
\|u_{tt}\|^2_{L^2(M)}+\|\nabla u_t\|^2_{L^2(M)}
&\leq&C[\|\varphi\|^2_{H^2(M)}+\|\psi\|^2_{H^1(M)}+\|f(0)\|^2_{L^2(M)}\nonumber\\&&+\int^t_0(\|f_t\|^2_{L^2(M)}
+\|\triangle u\|^2_{L^2(M)})d\tau]
\end{eqnarray}
Using (5.3)(5.5)(5.9)(5.10),one has
\begin{eqnarray}
\sum^2_{i=0}\|u\|_{W^{i,\infty}(I,H^{2-i}(M))}
&\leq&C[\|\varphi\|^2_{H^2(M)}+\|\psi\|^2_{H^1(M)}+\|f(0)\|^2_{L^2(M)}\nonumber\\&&+\int^T_0\|f_t\|^2_{L^2(M)}d\tau]
\end{eqnarray}
Here, $C$ is a constant depended only on manifold
$M,L,T,a,\alpha,\beta,\gamma$.

Then,use the Faedo-Galerkin method to give the existence of
(1.3).Let $\lambda_0,\lambda_1,\cdots$ be eginvalues of laplace
operator on Reimainnian manifold $M$,and $w_i$ be corresponded
eginvectors,namely $\triangle
w_i+\lambda_iw_i=0$,$i\in\{0,1,\cdots\}$.By Chapter 3 p87 of
[Schoen$\&$Yau],one knows that $w_i$ which are smooth functions on
$M$ are complete standard orthonormal basis of $L^2(M)$ and also
are complete orthogonal basis of $H^1(M)$,and $\lambda \geq
0$.Then we give approximate sequence
$u_n(x,t)=\sum^n_{i=0}\eta_i(t)w_i(x)$.In what follows,we denote
$(,)$ the inner product of $L^2(M)$.Consider \begin{eqnarray}
&&((u_n)_{tt},w_j)+(a(u_n)_t,w_j)-(\alpha\triangle
u_n+<\nabla\beta,\nabla u_n>+\gamma
u_n,w_j)\nonumber\\&=&(f,w_j),0\leq j\leq n
\end{eqnarray}
By the expenssion of $u_n$,
\begin{eqnarray}
&&\eta^{\prime\prime}_j(t)+\sum^n_{i=0}(aw_i,w_j)\eta^{\prime}_i+\sum^n_{i=0}(\alpha\lambda_iw_i-<\nabla\beta,\nabla
w_i>-\gamma w_i,w_j)\eta_i\nonumber\\&=&f_j:=(f,w_j),0\leq j\leq n
\end{eqnarray}
and on the condition,
\begin{eqnarray}
\eta_j(0)=\varphi_j:=(\varphi,w_j),\eta^{\prime}_j(0)=\psi_j:=(\psi,w_j)
\end{eqnarray}
This means
\begin{eqnarray}
u_n(0)=\sum^n_{i=0}\eta_i(0)w_i=\sum^n_{i=0}\varphi_iw_i,u^{\prime}_n(0)=\sum^n_{i=0}\psi_iw_i
\end{eqnarray}
So,let $A_{ij}(t)=(\alpha\lambda_iw_i-<\nabla\beta,\nabla
w_i>-\gamma w_i,w_j),B_{ij}(t)=(aw_i,w_j)$, then (5.13)(5.14) give a
linear ordinary differential system with initial data.By the theory
of ordinary differential equation.(5.13)(5.14) have a unique smooth
solution $(\eta_0,\cdots,\eta_n)$.

Obvious, $u_n(0)\longrightarrow\varphi$ in
${L^2(M)}$,$u^{\prime}_n(0)\longrightarrow\psi$ in $L^2(M)$.Then
using the closeness of $M$,
\begin{eqnarray}
\|\varphi\|^2_{H^1(M)}&=&(\varphi-\triangle\varphi,\varphi)=\sum^{\infty}_{i=0}(1+\lambda_i)|\varphi_i|^2<+\infty
\end{eqnarray}
This implies $u_n(0)\longrightarrow \varphi$ in $H^1(M)$.Similar
$u^{\prime}_n(0)\longrightarrow \psi$ in $H^1(M)$.And
\begin{eqnarray}
\|\varphi\|^2_{H^2(M)}&=&\sum^{\infty}_{i=0}(1+\lambda_i)|\varphi_i|^2+\int_M<-\triangle\nabla\varphi,\nabla\varphi>dV
\end{eqnarray}
By the Ricci identity,
\begin{eqnarray}
\int_M<-\triangle\nabla\varphi,\nabla\varphi>dV&=&\int_M[<-\nabla\triangle\varphi,\nabla\varphi>\nonumber\\&&-\sum^m_{i,j,p,q=1}g^{ij}g^{pq}R_{ip}\nabla_j\varphi\nabla_q\varphi]dV
\nonumber\\&\geq&(\triangle^2\varphi,\varphi)-C\|\varphi\|^2_{H^1(M)}\end{eqnarray}
Hence,$\sum^{\infty}_{i=0}(1+\lambda_i+\lambda^2_i)|\varphi|^2\leq
C\|\varphi\|^2_{H^2(M)}$.And,similar
\begin{eqnarray}
\|u_n(0)-\varphi\|^2_{H^2(M)}\leq\sum^{\infty}_{i=n+1}(1+\lambda_i+\lambda^2_i)|\varphi_i|^2+C\|u_n(0)-\varphi\|^2_{H^1(M)}
\longrightarrow 0\nonumber\\
\end{eqnarray}
when $n\longrightarrow\infty$.Then $u_n(0)\longrightarrow \varphi$
in $H^2(M)$. In (5.12),multiple $\eta^{\prime}_j$,and sum up $j$
from $0$ to $n$;multiple $\lambda_j\eta^{\prime}_j$, and sum up $j$
from $0$ to $n$;take derivation to $t$ then multiple
$\eta^{\prime\prime}_j$ and sum up $j$ from $0$ to $n$,similar to
(5.11) one obtains,
\begin{eqnarray}
\sum^2_{i=0}\|u_n\|_{W^{i,\infty}(I,H^{2-i}(M))}
&\leq&C[\|\varphi\|^2_{H^2(M)}+\|\psi\|^2_{H^1(M)}+\|f(0)\|^2_{L^2(M)}\nonumber\\&&+\int^T_0\|f_t\|^2_{L^2(M)}d\tau]
\end{eqnarray}
This means $u_n$ are uniformly bounded in
$\bigcap^2_{i=0}W^{i,\infty}(I,H^{2-i}(M))$.So exist some
$u\in\bigcap^2_{i=0}W^{i,\infty}(I,H^{2-i}(M))$, and
$u_n\longrightarrow u $ weakly star in $L^{\infty}(I,H^2(M));$
$(u_n)_t\longrightarrow u_t$
 weakly star in
$L^{\infty}(I,H^1(M))$; $(u_n)_{tt}\longrightarrow u_{tt}$ weakly
star in $L^{\infty}(I,L^2(M))$. This means
  $((u_n)_{tt},w_j)+(a(u_n)_t,w_j)-(\alpha\triangle u_n+<\nabla\beta,\nabla
u_n>+\gamma
u_n,w_j)\longrightarrow((u)_{tt},w_j)+(au_t,w_j)-(\alpha\triangle
u+<\nabla\beta,\nabla u>+\gamma u,w_j)$ weakly star in
$L^{\infty}[0,T]$.Hence $u$ is the solution of
(1.3).And by Lemma 2.7,$u$ satisfies (5.1).So one has\\

\noindent\textbf{Lemma 5.1}\quad\ Suppose $f\in
W^{1,\infty}(I,L^2(M)),a,\alpha,\beta,\gamma\in C^{\infty}(M\times
I),\varphi,\psi\in C^{\infty}(M)$,then there is a unique solution of
(1.3)(5.1) in $\bigcap^2_{i=0}W^{i,\infty}(I,H^{2-i}(M))$ and it has
prior estimate (5.11). Here, $(M,g)$ is a m-dimensional closed
Riemainnian manifold, $I=[0,T]$ is a interval of length $T$.\\

In what follows,we discuss the regularity of the solution of
(1.3).This comes from the following Lemma.\\

\noindent\textbf{Lemma 5.2}\quad\ The hypothesis of $M,I$ are same
to the above Lemma.And let $\forall N\in\mathbb{N},N\geq
2$,$a,\alpha,\beta,\gamma\in C^{\infty}(M\times I),\varphi,\psi\in
C^{\infty}(M),f\in\bigcap^{N-1}_{i=1}W^{i,\infty}(I,$
$H^{N-1-i}(M))$,then the solution of (1.3)(5.1)
$u\in\bigcap^N_{i=0}W^{i,\infty}(I,H^{N-i}(M))$,and $u$ satisfy the
following inequality:
\begin{eqnarray}
\sum^N_{i=0}\|\partial^i_tu\|^2_{H^{N-i}(M)}&\leq&C_N(\|\varphi\|^2_{H^N(M)}+\|\psi\|^2_{H^{N-1}(M)}
+\|f(0)\|^2_{H^{N-2}(M)}\nonumber\\&&+\int^t_0\sum^{N-1}_{i=1}\|\partial^i_tf\|^2_{H^{N-i-1}(M)}d\tau)
\end{eqnarray}
here,$C_N$ depends only on $N$ and manifold $M,L,T,a,\alpha,\beta,\gamma$.\\
\textbf{Proof}\quad\ We use induction to prove this Lemma.When
$N=2$,this is the result of Lemma 5.1.Then one suppose this is true
of case $N$,and consider the Lemma of case $N+1$.Take derivation to
$t$ in (1.3),and denote $v=u_t$,then
\begin{eqnarray}
&&v_{tt}+av_t-(\alpha\triangle v+<\nabla \beta,\nabla v>+\gamma
v)\nonumber\\&=&f_t-a_tv+\alpha_t\triangle(\varphi+\int^t_0vd\tau)+<\nabla
\beta_t,\nabla
(\varphi+\int^t_0vd\tau)>\nonumber\\&&+\gamma_t(\varphi+
\int^t_0vd\tau)
\end{eqnarray}
Introduce the following approximate sequence
$\{v_q\}^{\infty}_{q=0}$,
\begin{eqnarray}
&&(v_{q+1})_{tt}+a(v_{q+1})_t-(\alpha\triangle v_{q+1}+<\nabla
\beta,\nabla v_{q+1}>+\gamma
v_{q+1})\nonumber\\&=&f_t-a_tv_q+\alpha_t\triangle(\varphi+\int^t_0v_qd\tau)+<\nabla
\beta_t,\nabla
(\varphi+\int^t_0v_qd\tau)>\nonumber\\&&+\gamma_t(\varphi+
\int^t_0v_qd\tau)
\end{eqnarray}
on the condition that
  $v_0=0$, $v_{q+1}(x,0)=\psi(x)$, $v^{\prime}_{q+1}(x,0)=f(x,0)-a(x,0)\psi(x)+\alpha(x,0)\triangle\varphi(x)
$ $+<\nabla\beta(x,0),\nabla \varphi(x)>+\gamma(x,0)\varphi(x)$.Then
using (5.21)(5.23),one has
\begin{eqnarray}
&&\sum^N_{i=0}\|\partial^i_t(v_{q+1}-v_q)\|^2_{H^{N-i}(M)}\nonumber\\&\leq&C_N
\int^t_0\sum^{N-1}_{i=1}\|\partial^i_t[-a_t(v_q-v_{q-1})+\alpha_t\triangle\int^t_0(v_q-v_{q-1})d\tau^{\prime}
\nonumber\\&&+<\nabla\beta_t,\nabla\int^t_0(v_{q}-v_{q-1})d\tau^{\prime}>+\gamma_t\int^t_0(v_{q}-v_{q-1})d\tau^{\prime}]\|^2_{H^{N-i-1}(M)}d\tau\nonumber\\
\end{eqnarray}
Since,when $i\geq 1,q\geq 2$,
\begin{eqnarray}
&&\|\partial^i_t[\alpha_t\triangle\int^t_0(v_q-v_{q-1})d\tau^{\prime}]\|^2_{H^{N-i-1}(M)}\nonumber\\&\leq&
\sum^i_{j=0}\|\partial^j_t\int^t_0(v_q-v_{q-1})d\tau^{\prime}]\|^2_{H^{N-i+1}(M)}\nonumber\\
&\leq&\sum^i_{j=0}\underbrace{\int^t_0\cdots\int^t_0}_{i-j}\|\partial^i_t\int^t_0(v_q-v_{q-1})d\tau^{\prime}]\|^2_{H^{N-i+1}(M)}d\tau_1\cdots
d\tau_{i-j}\nonumber\\&\leq&C_N\int^t_0\|\partial^{i-1}_t(v_q-v_{q-1})\|^2_{H^{N-i+1}(M)}d\tau
\end{eqnarray}
Similar estimate is right to the other three terms of (5.24).Then using
(5.24)(5.25), one has,
\begin{eqnarray}
\sum^{N}_{i=0}\|\partial^i_t(v_{q+1}-v_q)\|^2_{H^{N-i}(M)}&\leq&C_N\int^t_0\sum^{N-1}_{i=0}\|\partial^i_t(v_{q}-v_{q-1})\|^2_{H^{N-i}(M)}
\end{eqnarray}
Denote
$h_{q+1}=\sum^{N-1}_{i=0}\|\partial^i_t(v_{q+1}-v_q)\|^2_{H^{N-i}(M)}$,the
above equation gives $h_{q+1}(t)$ $\leq C_N\int^t_0h_q(t)d\tau$,and
denote $C(T)=\|h_2\|_{L^{\infty}[0,T]}$.Then  by iteration,one can
obtain $h_{q+1}(t)\leq C(T)\frac{(C_Nt)^{q-1}}{(q-1)!}$.This
implies:
\begin{eqnarray}
\sum^{N}_{i=0}\|\partial^i_t(v_{q+1}-v_q)\|^2_{H^{N-i}(M)}&\leq&C(T)\frac{(C_Nt)^{q-1}}{(q-1)!}
\end{eqnarray}
Then,$\{v_q\}^{\infty}_{q=0}$ convergence to some $v$ in
$\bigcap^N_{i=0}W^{i,\infty}(I,H^{N-i}(M))$,and $v$ satisfy (5.22)
with the condition $v(0)=\psi,v^{\prime}(0)=f(0)-a(0)\psi+\alpha(0)\triangle
\varphi+<\nabla \beta(0),\nabla \varphi>+\gamma(0)\varphi$.Hence,$u\in \bigcap^N_{i=0}W^{i+1,\infty}(I,H^{N-i}(M))$.By
the $W^{k,p}$ estimate of elliptic equation,and (1.3), one has $u\in
\bigcap^{N+1}_{i=0}W^{i,\infty}(I,H^{N+1-i}(M))$, and $v$ satisfies,
\begin{eqnarray}
&&\sum^N_{i=0}\|\partial^i_tv\|^2_{H^{N-i}(M)}\nonumber\\&\leq&
C_N[\|\varphi\|^2_{H^{N+1}(M)}+\|\psi\|^2_{H^N(M)}+\|f(0)\|^2_{H^{N-1}(M)}
+\int^t_0\sum^{N-1}_{i=1}\|\partial^i_tf_t\|^2_{H^{N-i-1}(M)}d\tau\nonumber\\&&+\int^t_0
\sum^{N-1}_{i=1}\sum^i_{j=0}(\|\partial^j_tv\|^2_{H^{N-i-1}(M)}+\|\partial^j_t(\varphi+\int^t_0vd\tau^{\prime})\|^2_{H^{N-i+1}(M)})d\tau]\nonumber\\
&\leq&C_N[\|\varphi\|^2_{H^{N+1}(M)}+\|\psi\|^2_{H^N(M)}+\|f(0)\|^2_{H^{N-1}(M)}
+\int^t_0\sum^{N}_{i=2}\|\partial^i_tf\|^2_{H^{N-i}(M)}d\tau\nonumber\\&&+\int^t_0
\sum^{N-1}_{i=0}\|\partial^i_tv\|^2_{H^{N-i}(M)}d\tau]
\end{eqnarray}
Integrate it,and $v=u_t$, one has,
\begin{eqnarray}
\sum^N_{i=0}\|\partial^{i+1}_tu\|^2_{H^{N-i}(M)}&\leq&
C_N[\|\varphi\|^2_{H^{N+1}(M)}+\|\psi\|^2_{H^N(M)}+\|f(0)\|^2_{H^{N-1}(M)}
\nonumber\\&&+\int^t_0\sum^{N}_{i=1}\|\partial^i_tf\|^2_{H^{N-i}(M)}d\tau]
\end{eqnarray}
and,
\begin{eqnarray}
\|f\|^2_{H^{N-1}(M)}&\leq&\int^t_0\|f_t\|^2_{H^{N-1}(M)}d\tau+\|f(0)\|^2_{H^{N-1}(M)}
\end{eqnarray}
Using (1.3),the $W^{k,p}$ estimate of elliptic
equation,(5.29)(5.30),one can obtains the inequality (5.21) in
the case of $N+1$.This complete the proof.\\

\noindent\textbf{Proposition 5.3}\quad\ The hypothesis of $M,I$ as
 above.Suppose $f,a,\alpha,\beta,\gamma\in C^{\infty}(M\times
 I),\varphi,\psi\in C^{\infty}(M)$.Then (1.3)(5.1) have a unique smooth solution and it satisfies
 (5.21).\\
 \textbf{Proof}\quad\ It is obvious from Lemma 5.1 and 5.2.
\section{Solution of Equation in Some Cases}

\quad\ The idea of this section can be found in
[H\"{o}mander],[John],[Klainerman]. Suppose $(M,g)$ be m-dimensional
closed Riemainnian manifold,$I=[0,T]$ is a interval of length $T$.By
Proposition 5.3,one can introduce the Picard iterate sequence
$\{u_n\}^{+\infty}_{n=1}$ as follows,
\begin{eqnarray}
m(u_{n+1})_{tt}-e^{-2(m+1)u_n}\triangle u_{n+1}&=&F(u_{n+1},u_n)
\end{eqnarray}
with the condition
$u_{n+1}(x,0)=\varphi(x),(u_{n+1})_t(x,0)=\psi(x),\varphi,\psi\in
C^{\infty}(M),u_1$ $=0$.Here, $F$ is according to (4.1).Let
$s>\frac{m}{2}+1$,and we suppose that exist $\tilde{T}\leq T$ such
that in $[0,\tilde{T}]$,exist positive constant $D$ satisfying,to
any $n$,
\begin{eqnarray}
\|u_n\|_{H^{s+1}(M)}+\|(u_n)_t\|_{H^s(M)}&\leq&D
\end{eqnarray}
Let
\begin{eqnarray}
G(u,v)&=&[1-e^{2(m+1)(u-v)}](R_g+2\triangle
u)-(m^2+m+2)[<\nabla(u-v),\nabla
u>\nonumber\\&&+(1-e^{2(m+1)(u-v)})<\nabla v,\nabla
u>]-m(3m+1)e^{2(m+1)u}u_t(u-v)_t\nonumber\\
\end{eqnarray}
Then,
\begin{eqnarray}
&&m(u_{n+1}-u_n)_{tt}-e^{-2(m+1)u_n}\triangle(u_{n+1}-u_n)\nonumber\\&=&
F(u_{n+1},u_n)-F(u_n,u_{n-1})+(e^{-2(m+1)u_n}-e^{-2(m+1)u_{n-1}})\triangle
u_n\nonumber\\&=&
\frac{1}{2}\{-e^{-2mu_n}(1-e^{2m(u_n-u_{n-1})})\tilde{R}+e^{-2(m+1)u_n}G(u_n,u_{n-1})\nonumber\\&&
-(m^2+m+2)e^{-2(m+1)u_n}<\nabla
u_n,\nabla(u_{n+1}-u_n)>\nonumber\\&&-m(3m+1)(u_n)_t(u_{n+1}-u_n)_t\}
\end{eqnarray}
with the condition that $(u_{n+1}-u_n)(x,0)=(u_{n+1}-u_n)_t(x,0)=0
$.So,we replace $u,v,\tilde{R},R_g$ by
$u_{n+1}-u_n,u_n,(1-e^{2m(u_n-u_{n-1})})\tilde{R},G(u_n,u_{n-1})$
respectively in (4.1),using (4.9),one has,
\begin{eqnarray}
A_s&=&\|(1-e^{2m(u_n-u_{n-1})})\tilde{R}\|_{H^{s-1}(M)}+\|G(u_n,u_{n-1})\|_{H^{s-1}(M)}
\end{eqnarray}
and using (6.2),and Corollary 2.4,and Sobolev embedding theorems,
\begin{eqnarray}
&&\|(1-e^{2m(u_n-u_{n-1})})\tilde{R}\|_{H^{s-1}(M)}\nonumber\\&\leq
&C\|1-e^{2m(u_n-u_{n-1})}\|_{H^{s-1}(M)}\|\tilde{R}\|_{H^{s-1}(M)}\nonumber\\
&\leq&Ce^{C\|u_n-u_{n-1}\|_{H^{s-1}(M)}}\|u_n-u_{n-1}\|_{H^{s-1}(M)}\|\tilde{R}\|_{H^{s-1}(M)}\nonumber\\
&\leq&Ce^{CD}\|u_n-u_{n-1}\|_{H^{s-1}(M)}\|\tilde{R}\|_{H^{s-1}(M)}
\end{eqnarray}
and similar,
\begin{eqnarray}
&&\|G(u_n,u_{n-1})\|_{H^{s-1}(M)}\nonumber\\&\leq&Ce^{CD}\|u_n-u_{n-1}\|_{H^{s-1}(M)}(\|R_g\|_{H^{s-1}(M)}+D)
+CD\|u_n-u_{n-1}\|_{H^{s}(M)}\nonumber\\&&+Ce^{CD}\|u_n-u_{n-1}\|_{H^{s-1}(M)}D^2+CD(1+e^{CD})\|
(u_n-u_{n-1})_t\|_{H^{s-1}(M)}\nonumber\\
\end{eqnarray}
Then,By (6.5)-(6.7),
\begin{eqnarray}
A_s&\leq&Ce^{CD}(D+\|\tilde{R}\|_{H^{s-1}(M)}+\|R_g\|_{H^{s-1}(M)})\nonumber\\&&\times(\|u_n-u_{n-1}\|_{H^s(M)}
+\|(u_n-u_{n-1})_t\|_{H^{s-1}(M)})
\end{eqnarray}
Then,using (4.33) and (6.2),
\begin{eqnarray}
\frac{d
(E_s(t))^{\frac{1}{2}}}{dt}&\leq&Ce^{CD}(A_s+(E_s(t))^{\frac{1}{2}})
\end{eqnarray}
here,do correspond replacement which is stated above in
$E_s(t)$.Then,integrate (6.9),and using (6.8),one obtains
\begin{eqnarray}
&&\|u_n-u_{n-1}\|_{H^s(M)}+\|(u_n-u_{n-1})_t\|_{H^{s-1}(M)}\nonumber\\
&\leq&C_4e^{C_5e^{C_6D}T}(D+\|\tilde{R}\|_{C^0(I,H^{s-1}(M))}+\|R_g\|_{H^{s-1}(M)})\nonumber\\&&\times
\int^t_0(\|u_n-u_{n-1}\|_{H^s(M)}+\|(u_n-u_{n-1})_t\|_{H^{s-1}(M)})
\end{eqnarray}
here, $C_4,C_5,C_6$ only relates to manifold $M$ and $s$.\\

\noindent\textbf{Proposition 6.1}\quad\ Suppose $(M,g)$ be
m-dimensional closed Riemainnian manifold, $I=[0,T]$ be a interval
of length $T$, and $\tilde{R}\in C^{\infty}(M\times I)$.To any
integer $k\geq 2$,exist $0< t_0\leq T$,such that equation (1.2) has
a solution $u\in C^k(M\times [0,t_0])$,and it satisfies
$u(x,0)=\varphi(x),u_t(x,0)=\psi(x)$.Here, $\varphi(x),\psi(x)\in
C^{\infty}(M)$ are two prescribed functions.And $t_0$ relates to
manifold $M,k,\tilde{R},\varphi,\psi$.\\
\textbf{Proof}\quad\ Take $s=[\frac{m}{2}+k+1]$ (i.e.the biggest
integer not bigger than $\frac{m}{2}+k+1$ ),then
$H^s(M)\hookrightarrow C^k(M)$. By Proposition 4.1(i),to any given
$\varphi,\psi$,take $D$ to satisfies
\begin{eqnarray}
D&>&2\sqrt{2}[1+(\sum^{s+1}_{l=1}\frac{1}{2}\int_M[|\nabla^l\varphi|^2+me^{2(m+1)\varphi}|\nabla^{l-1}\psi|^2+e^{-2(m+1)\varphi}|\nabla^l\varphi|^2\nonumber\\&&
+m|\nabla^{l-1}\psi|^2]dV+\frac{1}{2}\int_M|\varphi|^2dV)^{\frac{1}{2}}]
\end{eqnarray}
then exist $t_1$,in $[0,t_1]$,(4.34) implies (4.35).To
$\{u_n\}^{\infty}_{n=1}$ which is constructed in (6.1), obvious
$D\geq \|u_1\|_{H^{s+1}(M)}+\|(u_1)_t\|_{H^{s}(M)}=0$,so by
iteration,one has (6.2).Then take $t_0\leq t_1$  small enough such
that,
\begin{eqnarray}
C_4e^{C_5e^{C_6D}T}(D+\|\tilde{R}\|_{H^{s-1}(M)}+\|R_g\|_{H^{s-1}(M)})t_0&<&1
\end{eqnarray}
$C_4,C_5,C_6$ are constants in (6.10).Then using (6.10),
\begin{eqnarray}
&&\|u_n-u_{n-1}\|_{C^0([0,t_0],H^s(M))}+\|(u_n-u_{n-1})_t\|_{C^0([0,t_0],H^{s-1}(M))}
\nonumber\\&\leq&C_4e^{C_5e^{C_6D}T}(D+\|\tilde{R}\|_{C^0(I,H^{s-1}(M))}+\|R_g\|_{H^{s-1}(M)})t_0\nonumber\\&&\times
(\|u_n-u_{n-1}\|_{C^0([0,t_0],H^s(M))}+\|(u_n-u_{n-1})_t\|_{C^0([0,t_0],H^{s-1}(M))})\nonumber\\
\end{eqnarray}
This implies $\{u_n\}^{+\infty}_{n=1}$ and
$\{(u_n)_t\}^{+\infty}_{n=1}$ are convergence to some $u,u_t$ in
$C^0([0,t_0],$ $H^s(M))$,$C^0([0,t_0],H^{s-1}(M))$ respectively.
Obvious $u_n\longrightarrow u$ in $ C^0([0,t_0],$ $C^k(M))\bigcap$ $
C^1([0,t_0],C^{k-1}(M))$.By equation (6.1)
$(u_n)_{tt}\longrightarrow u_{tt}$ in $C^0([0,t_0],$
$C^{k-2}(M))$.Then $u\in C^2([0,t_0],C^{k-2}(M))$ and $u$ satisfies
equation (1.2).Take derivation to $t$ of first order,one can finds
$u\in
C^3([0,t_0],C^{k-3}(M))$.Then iterate this step shows the result holds.\\

\noindent\textbf{Proposition 6.2}\quad\  Suppose $(M,g)$ be
m-dimensional closed Riemainnian manifold which is scalar curvature
flat $R_g=0$.And $I=[0,T]$ be a interval of length $T$, and
$\tilde{R}\in C^{\infty}(M\times I)$.To any integer $k\geq 2$,exist
$\varepsilon
>0$,such that equation (1.2) has a solution $u\in C^k(M\times
I)$,and it satisfies $u(x,0)=0,u_t(x,0)=0$ when
$\|\tilde{R}\|_{H^{[\frac{m}{2}+k]}(M)}\leq \varepsilon$.Here,
$\varepsilon$ only relates to
manifold $M,k,I$.\\
\textbf{Proof}\quad\ Similar to Proposition 6.1,take
$s=[\frac{m}{2}+k+1]$,then $H^s(M)\hookrightarrow C^k(M)$.One can
take $D,\varepsilon$ small enough such that,
\begin{eqnarray}
C_4e^{C_5e^{C_6D}T}(D+\varepsilon)&<&1
\end{eqnarray}
$C_4,C_5,C_6$ are constants in (6.10) and Proposition 4.2 holds.
Then (4.44) implies (4.45) and the construction of
$\{u_n\}^{+\infty}_{n=1}$ in (6.1),one finds (6.2).This implies
(6.10),namely,
\begin{eqnarray}
&&\|u_n-u_{n-1}\|_{C^0(I,H^s(M))}+\|(u_n-u_{n-1})_t\|_{C^0(I,H^{s-1}(M))}
\nonumber\\&\leq&C_4e^{C_5e^{C_6D}T}(D+\varepsilon)
(\|u_n-u_{n-1}\|_{C^0(I,H^s(M))}\nonumber\\&&+\|(u_n-u_{n-1})_t\|_{C^0(I,H^{s-1}(M))})
\end{eqnarray}
Then, similar argument of Proposition 6.1 gives the result by
(6.14)(6.15).\\

 At the end we discuss the uniqueness of solution in the above
two Proposition.This is from the next result. \\

\noindent\textbf{Proposition 6.3}\quad\ Suppose $(M,g)$ be
m-dimensional closed Riemainnian manifold, $I=[0,T]$ be a interval
of length $T$, and $\tilde{R}\in C^{\infty}(M\times I)$.Let $0\leq
\tilde{T}\leq T$.If equation (1.2) has a solution $u\in
C^{k+2}(M\times [0,\tilde{T}])$,and $u$ satisfies
$u(x,0)=\varphi(x),u_t(x,0)=\psi(x)$,then the solution is unique.
Here,$\varphi(x),\psi(x)\in C^{\infty}(M)$ are two prescribed
functions,and $k>\frac{m}{2}+1,k\in\mathbb{Z}$.\\
\textbf{Proof}\quad\ Suppose $u,v\in C^{k+1}(M\times [0,\tilde{T}])$
are two different solution of (1.2) with the same initial
data.So,exist $t_{min}$,$0<t_{min}<\tilde{T}$ and in
$[0,t_{min}],u=v$,but $u\neq v$ when $t>t_{min}$.Let
$\tilde{u}(x,t)=u(x,t+t_{min}),\tilde{v}(x,t)=v(x,t+t_{min})$,then
$\tilde{u}(x,0)=\tilde{v}(x,0),\tilde{u}_t(x,0)=\tilde{v}_t(x,0)$.And
$\tilde{u},\tilde{v}$ also satisfy (1.2),differ the corresponded
equation,one has
\begin{eqnarray}
&&(\tilde{u}-\tilde{v})_{tt}-e^{-2(m+1)\tilde{u}}\triangle(\tilde{u}-\tilde{v})\nonumber\\&=&
\frac{1}{2}\{-e^{-2m\tilde{u}}(1-e^{2m(\tilde{u}-\tilde{v})})\tilde{R}+e^{-2(m+1)\tilde{u}}H(\tilde{u},\tilde{v})\nonumber\\&&
-(m^2+m+2)e^{-2(m+1)\tilde{u}}<\nabla
\tilde{u},\nabla(\tilde{u}-\tilde{v})>-m(3m+1)\tilde{u}_t(\tilde{u}-\tilde{v})_t\}\nonumber\\
\end{eqnarray}
where,
\begin{eqnarray}
H(\tilde{u},\tilde{v})&=&[1-e^{2(m+1)(\tilde{u}-\tilde{v})}](R_g+2\triangle
\tilde{v})-(m^2+m+2)[<\nabla(\tilde{u}-\tilde{v}),\nabla
\tilde{v}>\nonumber\\&&+(1-e^{2(m+1)(\tilde{u}-\tilde{v})})<\nabla
\tilde{v},\nabla
\tilde{v}>]-m(3m+1)e^{2(m+1)\tilde{u}}\tilde{v}_t(\tilde{u}-\tilde{v})_t\nonumber\\
\end{eqnarray}
Similar to (6.2)-(6.10) and using Proposition 4.1(ii),one obtains,
\begin{eqnarray}
&&\|\tilde{u}-\tilde{v}\|_{H^k(M)}+\|(\tilde{u}-\tilde{v})_t\|_{H^{k-1}(M)}
\nonumber\\&\leq&C
\int^t_0(\|\tilde{u}-\tilde{v}\|_{H^k(M)}+\|(\tilde{u}-\tilde{v})_t\|_{H^{k-1}(M)})d\tau
\end{eqnarray}
holds in $[0,t_0]$,and $t_0$ is small enough.Constant $C$ and $t_0$
relates to manifold
$M,k,\tilde{R},\tilde{u}(0),\tilde{u}_t(0)$.(6.18) implies $u=v$
when $t_0$ more small such that $Ct_0<1$.This is a
contradiction.So,$t_{min}=\tilde{T}$.\\

Proposition 6.1-6.3 imply Theorem 1.1 and 1.2.\\

At last,we suggest some questions which maybe is nature to ask:\\
(1)Replace $I$ by $S^1$ (i.e. the unit cycle),what is the sufficient
and necessary condition on which equation (1.2) exists solution ?\\
(2)On what condition,equation (1.2) exists solution on $M\times
\mathbb{R}^1$ ?\\
(3)How to solve this question of the higher dimension version,namely
solve equation (3.31) ?\\
(4)Can this question be generalized to a fiber bundle over a base
manifold $M$ ? What is the corresponded equation ?\\

\begin{center}
{\textbf{References}}
\end{center}
\noindent[Chern,Chen$\&$Lam] S.S.Chern,W.H.Chen,K.S.Lam,Lectures on
Differential Geometry,Series on University Mathematics Vol 1,World
Scientific,1998.

\noindent[Gilbarg$\&$Trudinger] D.Gilbarg,N.S.Trudinger, Elliptic
Partial Differential Equations of Second Order,Classics in
Mathematics,Springer,2000.

\noindent [H\"{o}rmander] L.H\"{o}rmander, Lectures on Nonlinear
Hyperbolic Differential Equation,Math\'{e}matiques $\&$ Applications
26,Springer-Verlag,1997.

\noindent [Hughes,Kato$\&$Marsden] T.J.R.Hughes,T.Kato,and
J.E.Marsden,Well-Posed Quasi-Linear Second -Order Hyperbolic Systems
with Application to Nonlinear Elastodynamics and General
Relativity,Arch.Rat.Mech.Anal.63(1977),pp 273-294.

\noindent [John] F.John,Decayed Singularity Formation in Solution of
Nonlinear Wave Equations in Higher Dimension.Comm.Pure
Appl.Math.29(1976),pp 649-681.

\noindent [Klainerman] S.Klainerman,Global Existence for Nonlinear
Wave Equations. Comm.Pure Appl.Math.1980,Vol \textrm{XXXIII},pp
43-401.

\noindent[Nirenberg] L.Nirenberg, On Elliptic Partial Differential
Equation, Ann.Scu. Norm.Sup.Pisa 13(3)(1959),115-162.

\noindent [Smith \& Tataru] H.F.Smith and D.Tataru,Sharp Local
Well-Posedness Results for the Nonlinear Wave Equation.Ann.of
Math.163(2005),pp291-366.

\noindent [Xin] J.Xin,Doctor Thesis.

\noindent [Schoen$\&$Yau] R.Schoen,S-T.Yau, Lectures on Differential
Geometry, Conference Proceedings and Lecture Notes in Geometry and
Topology,Volume 1,International Press,1994.

\noindent[Zheng] S.Zheng, Nonlinear Evolution Equations, Chapman
$\&$ Hall/CRC Monographs and Surveys in Pure and Applied Mathematics
133, Chapman $\&$  Hall/CRC,2004.





\end{document}